\documentclass[a4paper, 11pt]{article}
\usepackage{amsmath} \usepackage{amsfonts} \usepackage{amssymb}\usepackage{latexsym}
\usepackage[english]{babel}
\usepackage{amscd}
\usepackage[dvips,final]{graphics}
\usepackage{locengli,math}
\usepackage{graphicx,pst-plot,psfig,pstricks,pst-node,pst-text,pst-tree,pst-char,pst-coil}
\begin{document}
\begin{center}
\textbf{\LARGE{\textsf{On some remarkable operads constructed from Baxter operators}}}
\footnote{
{\it{2000 Mathematics Subject Classification: 17A30, 18D50.}}
{\it{Key words and phrases: dendriform dialgebras, quadri-algebras, octo-algebras, dipterous algebras, coassociative codialgebras, coassociative cotrialgebras, Leibniz algebras, Poisson dialgebras, Baxter operators.}} 
}
\vskip1cm
\parbox[t]{14cm}{\large{
Philippe {\sc Leroux}}\\
\vskip4mm
{\footnotesize
\baselineskip=5mm
Institut de Recherche
Math\'ematique, Universit\'e de Rennes I and CNRS UMR 6625\\
Campus de Beaulieu, 35042 Rennes Cedex, France, pleroux@univ-rennes1.fr}}
\end{center}
\vskip1cm
\baselineskip=5mm
\noindent
\begin{center}
\begin{large}
\textbf{ 20/10/03}
\end{large}
\end{center}
{\bf Abstract:}
Language theory, symbolic dynamics, modelisation of viral insertion into the genetic code of a host cell motivate the introduction of new types of bialgebras, whose coalgebra parts are not necessarily coassociative. Similarly, certain algebraic descriptions of combinatorial objects such as weighted directed graphs through a coalgebraic formulation often require at least two coalgebras whose coproducts verify 
entanglement conditions. All these structures generate special types of associative algebras (or binary quadratic and non-symmetric operads) obtained from commuting Baxter operators. One of the aim of this article is to study what type of associative algebras appear when such or such coalgebraic structures are used to describe for instance combinatorial objects. Examples from symbolic dynamics, weighted directed graphs and matrix theory are given.
 
\begin{scriptsize}
\tableofcontents
\end{scriptsize}
\section{Introduction and notation}
In the sequel $k$ is a field of characteristic zero.
Let $(X, \bullet)$ be a $k$-algebra and $(\bullet_i)_{1 \leq i \leq N}: X^{\otimes 2} \xrightarrow{} X$ be a family of binary operations on $X$. 
The notation $\bullet \longrightarrow \sum_i \bullet_i$ will mean $x \bullet y =  \sum_i x \bullet_i y$, for all $x,y \in X$. We say that
the operation $\bullet$ {\it{splits}} into the $N$ operations $\bullet_1, \ldots, \bullet_N,$ or that the operation $\bullet$ is a {\it{cluster of $N$ (binary) operations}}.

\noindent
This paper is a continuation of \cite{AguiarLoday, Lertribax} on $\epsilon$-bialgebras and commuting Baxter operators on the one hand and of \cite{dipt,codialg1,Coa} on weighted directed graphs on the other hand. A $\epsilon$-bialgebra
$(A,\mu,\Delta)$ is an associative algebra $(A,\mu)$ together with a coassociative coalgebra $(A,\Delta$) such that, for all $a,b \in A$, $\Delta(ab) := \Delta(a)b + a \Delta(b)  .$ Such bialgebras
appeared for the first time in the work of Joni and Rota \cite{Rota}, see also \cite{Aguiar, AguiarLoday, Lertribax, Voi}.
However,
other structures, corresponding to the equality $\Delta(ab) := \Delta(a)b$ or $\Delta(ab) := a \Delta(b)$ can be met in algebraic descriptions of some substitutions in symbolic dynamics (see Section 1) or in the work of Khovanov \cite{Khovanov} on 2-dimensional topological quantum field theory. They are called left or right $\epsilon$-bialgebras and are defined as follows.
\begin{defi}{[$\epsilon[R]$-, $\epsilon[L]$-Bialgebra]} A {\it{$\epsilon[R]$-bialgebra}} (resp. {\it{$\epsilon[L]$-bialgebra}})        
$(A,\mu,\Delta)$ is an associative algebra $(A,\mu)$ together with a coassociative coalgebra $(A,\Delta$) such that, for all $a,b \in A$, $\Delta(ab) := a \Delta(b)$, (resp. $\Delta(ab) := \Delta(a) b$). Similarly, a {\it{$\epsilon'[R]$-bialgebra}} (resp. {\it{$\epsilon'[L]$-bialgebra}}) $(A,\mu,\Delta)$ is an associative algebra $(A,\mu)$ together with a coalgebra $(A,\Delta$) not supposed to be coassociative such that, for all $a,b \in A$, $\Delta(ab) := a \Delta(b)$, (resp. $\Delta(ab) := \Delta(a) b$). 
The linear map $\Delta: A \xrightarrow{} A^{\otimes 2}$ is said to be a {\it{co-operation}}.
\end{defi}
To any $\epsilon$-bialgebras, commuting pairs of Baxter operators can be 
introduced \cite{AguiarLoday}. 
Let $(A, \ \mu)$ be an associative algebra. Recall that a Baxter operator is a linear map $\zeta: A \xrightarrow{} A$
such that, for all $x,y \in A,$ 
$$ \zeta(x) \zeta(y)=\zeta(\zeta(x) y + x  \zeta(y)).$$
Baxter opeartors on associative algebras appear originally in a work of G. Baxter \cite{Baxter} and the importance of such a map was stressed by G.-C. Rota in \cite{Rota1}. See also \cite{Aguiar, AguiarLoday, KEF, Lertribax} for more details.

\noindent
Similarly, the notion of right Baxter operators (resp. left Baxter operators), i.e., linear maps $\zeta: A \xrightarrow{} A$ verifying
$\zeta(x)\zeta(y)=\zeta(\zeta(x) y))$ (resp. $\zeta(x) \zeta(y)=\zeta(x \zeta(y))$) on associative algebras can be introduced.
Such operators appear in the study of $\epsilon'[R]$-bialgebras and 
their left versions~\footnote{In the sequel of this introduction, we will focus mainly on right versions.}. 
The novelty lies in $\epsilon'[R]$-bialgebras where co-operations are not supposed to be coassociative since right Baxter operators can be introduced as well as associative structures, called $L$-anti-dipterous algebras on $\textsf{End}(A)$, the $k$-algebra of linear maps from $A$ to $A$. The case of $\epsilon[R]$-bialgebras, where co-operations are supposed to be coassociative, produce left Baxter operators which commute with previous right ones. Other types of associative algebras, presenting symmetries, will appear.
This is studied in Section 2. Examples are given. We also recall
the definition of quadri-algebras \cite{AguiarLoday} and give other examples of $\epsilon$-bialgebras. In an example already described in \cite{Lertribax}, we relate $\epsilon$-bialgebra structures to a toy-model of virus attacks inside a host cell.

To describe some mathematical objects, such as weighted directed graphs \cite{Coa} or how graphs can be covered by coassociative structures, \cite{dipt, codialg1},  is possible only with the help of several co-operations. That is why
in Section 3 and 4, we explore what occurs if associative algebras have two co-operations, related by an entanglement condition, which well behave with the underlying product. Such descriptions also produce commuting Baxter operators. Many cases are undertaken. This leads to the constructions of $n$-hypercubic algebras, $n$-circular algebras, both related to the notion of Hochschild 2-cocycles. 

Up to now, we have worked with two commuting Baxter operators and their right or left versions. What happens if three Baxter operators commute pairwise? In section 5, we introduce particular $k$-vector spaces with two coproducts, called $L$-circular $\epsilon$-bialgebras. Via such structures, three pairwise commuting Baxter operators are created. In the 
case of $\epsilon$-bialgebras, two commuting Baxter operators were created giving birth to quadri-algebras \cite{AguiarLoday}, i.e., $k$-vector spaces $Q$ equipped with four binary operations $\nwarrow, \ \nearrow, \ \swarrow, \ \searrow \ : \ Q^{\otimes 2} \xrightarrow{} Q$ and obeying 9 compatibility axioms. One of the interests of such algebras is to produce an associative algebra $(Q, \star)$ such that $\star$ is the following associative cluster $\star  \xrightarrow{} \ \nwarrow + \nearrow + \swarrow + \searrow$. The second and main interest is to produce two dendriform dialgebras by combining these four operations, the two dendriform dialgebras $(Q, \prec, \succ)$ and $(Q, \wedge, \vee)$ being such that $\star \xrightarrow{} \ \prec + \succ$ and $\star \xrightarrow{} \  \wedge + \vee$. In the case of three pairwise commuting Baxter operators, the relevant notion is called octo-algebras. Octo-algebras are associative algebras whose associative product is a cluster of eight operations. By combining them, three quadri-algebras (or six dendriform dialgebras) sharing the same associative product can be produced. Examples are given, one of them coming from discrete mathematics.
This work ends with Section~6 by showing how to produce, from right Baxter operators, Leibniz and Poisson dialgebras, notion introduced in \cite{Loday, LodayRonco}. Examples from coassociative coverings of some directed graphs are also given.

\noindent
The background to make this article as self-contained as possible is now presented.
\begin{defi}{[Dendriform dialgebra \cite{Loday}]}
A {\it{dendriform dialgebra}} is a $k$-vector space $E$ equipped with two binary operations:
$\prec, \ \succ: E^{\otimes 2} \xrightarrow{} E$, satisfying the following relations for all $x,y \in E$:
$$(x \prec y )\prec z = x \prec(y \star z), \ \ \
(x \succ y )\prec z = x \succ(y \prec z), \ \ \
(x \star y )\succ z = x \succ(y \succ z), \ \ \ $$
where by definition $x \star y :=x  \prec y +x \succ y$, for all $x,y \in E$. The dendriform dialgebra $(E,\star)$ is then an associative algebra such that $\star \longrightarrow \ \prec + \succ$. Observe that these axioms are globally invariant under the transformation $x \prec^{op} y := y \succ x$, $x \succ^{op} y := y \prec x$. A dendriform dialgebra is said to be {\it{commutative}} if $x \prec y := y \succ x$. Such algebras are also called {\it{Zinbiel algebras}} \cite{Loday}. The free dendriform dialgebra is related to planar (rooted) binary trees.
\end{defi}
\begin{defi}{[$L$-Algebra]}
A {\it{$L$-algebra}} $(A, \succ \ \dashrightarrow \ \prec)$ is a $k$-vector space
equipped with two binary operations, $\succ,  \ \prec: A^{\otimes 2} \xrightarrow{} A$ {\it{linked}} for all $x,y,z \in A$ by,
$$ (L.) \ \ \ (x \succ y) \prec z=x \succ (y \prec z) \equiv \ \succ \ \dashrightarrow \ \prec.$$
A $L$-algebra is said to be {\it{associative}} if both $\prec$ and $\succ$ are associative.
The free $L$-algebra is related to planar (rooted) binary trees \cite{Pir}. The importance of $L$-coalgebras can be found in \cite{Coa} and in the notion of coverings or tilings of directed graphs by coassociative structures \cite{dipt,codialg1}. By the notation $\prec \dashrightarrow \succ$, we mean $(x \prec y) \succ z=x \prec (y \succ z).$ Similarly, for the coversion,  $\Delta_{\prec} \dashrightarrow \Delta_{\succ}$ means
$(\Delta_{\prec} \otimes id) \Delta_{\succ}:=(id \otimes \Delta_{\succ}) \Delta_{\prec} $. Such an equality will be called an {\it{entanglement}}. 
\end{defi}
\begin{defi}{[$L$-Dipterous algebra]}
\label{LAD}
A {\it{dipterous algebra}} $(A, \star, \succ)$ is an associative algebra
$(A, \star)$ equipped with a left module on itself, i.e., for all $x,y,z \in A$,
$$ (As.) \ \ \ (x \star y) \star z=x \star (y \star z), \ \ \
(Dip.) \ \ \ (x \star y) \succ z=x \succ (y \succ z). $$
The free dipterous algebra is related to two copies of planar rooted trees \cite{LR1}.
Similarly, a {\it{$L$-dipterous algebra}} $(A, \star, \succ)$ is a dipterous algebra verifying an extra condition:
$$(L.)  \ \ \ (x \succ y) \star z = x \succ (y \star z),$$
for all $x,y \in A$.
A {\it{$L$-anti-dipterous}} algebra is a $k$-vector space $A$ equipped with
two binary operations $\star, \ \prec: A^{\otimes 2} \xrightarrow{} A,$ obeying, 
$$ (As.) \ \ (x \star y) \star z=x \star (y \star z), \ \ \ (A.Dip.) \ \ (x \prec y) \prec z=x \prec (y \star z), \ \
(L.) \ \ (x \star y)\prec  z = x \star (y \prec z).$$
\end{defi}
In terms of category, we get,
$$\textsf{L-Dip.}  \longrightarrow \textsf{Dip.} \longrightarrow As.
 \ \ \ \textrm{and} \ \ \ \textsf{L-A.Dip.}  \longrightarrow \textsf{A.Dip.} \longrightarrow As.$$
\noindent
All these $k$-vector spaces are said to be algebras of a certain type $P$. To better formulate this sentence, we need the operadic language, see for instance \cite{Fresse,GK,Lodayscd}. Let us briefly expose what this is.  
Let $P$ be a type of algebras and $P(V)$ be the free $P$-algebra on the $k$-vector space $V$, i.e., by definition,
a $k$-algebra of type $P$ equipped with a map $i: \ V \mapsto P(V)$ which satisfies the following universal property:
for any linear map $f: V \xrightarrow{} A$, where $A$ is a $k$-algebra of type $P$, there exists a unique morphism of algebra of type $P$, $\bar{f}: P(V) \xrightarrow{} A$ such that $\bar{f} \circ i=f$.
Suppose $P(V) := \oplus_{n \geq 1} P(n) \otimes _{S_n} V^{\otimes n},$ where $P(n)$ are right $S_n$-modules. We view $P$ as an endofunctor on the category of $k$-vector spaces. The structure of the free $P$-algebra of $P(V)$ induces a natural transformation $\pi: P \circ P \xrightarrow{} P$ as well as $u: Id \xrightarrow{} P$ verifying usual associativity and unitarity axioms. An {\it{algebraic operad}} is then a triple $(P, \ \pi, \ u)$. A {\it{$P$-algebra}} or algebra of type $P$ is then a $k$-vector space $A$ together with a linear map $\pi_A: P(A) \xrightarrow{} A$ such that $\pi_A \circ \pi(A) = \pi_A \circ P(\pi_A)$ and $\pi_A \circ u(A) = Id_A.$
The $k$-vector space $P(n)$ is the space of $n$-ary operations for $P$-algebras. In the sequel we suppose there is, up to homotheties, a unique $1$-ary operation, the identity, i.e.,  $P(1):= k.Id$ and that all possible operations are generated by composition from $P(2)$, i.e., by binary operations. The operad is said to be {\it{binary}}. It is said to be {\it{quadratic}} if all the relations between operations are consequences of relations described exclusively with the help of monomials with two operations. An operad is said to be {\it{non-symmetric}} if, in the relations, the variables $x,y,z$ appear in the same order. In this case, the $k$-vector space $P(n)$ can be written like $P(n):= P'(n) \otimes k[S_n]$, where $P'(n)$ is also a $k$-vector space.
The free $P$-algebra is then entirely induced by the free $P$-algebra on one generator $P(k):= \oplus_{n \geq 1} \ P'(n)$. The generating function of the operad $P$ is given by:
$$ f^{P}(x):= \sum \ (-1)^n \frac{\textrm{dim} \ P(n)}{n!} x^n := \sum \ (-1)^n \textrm{dim} \ P'(n) x^n.$$
In this paper, operads will be binary, quadratic and non-symmetric. Their sequences $(\textrm{dim} \ P'(n))_{n \geq 1}$, at least their starts will be indicated.

\begin{defi}{[Baxter operators]}
Let $(P, \pi, u)$ be a binary, quadratic and non-symmetric operad. A {\it{Baxter operator}}, (resp. {\it{right Baxter operator}}), (resp. {\it{left Baxter operator}}) on the $k$-algebra $A$ of type $P$, is a linear map $\zeta: A \xrightarrow{} A$, such that, for all $x,y \in A$ and for all $\circ \in P(2)$,
$$ \zeta(x)\circ \zeta(y)=\zeta(\zeta(x)\circ y + x \circ \zeta(y)),$$
(resp. $\zeta(x)\circ \zeta(y)=\zeta(\zeta(x)\circ y)$), (resp. $\zeta(x)\circ \zeta(y)=\zeta(x\circ \zeta(y))$.
\end{defi}
If $(A, \ \mu)$ is an associative $k$-algebra, then 
the following operations $a (x \otimes y) := ax \otimes y$ and $ (x \otimes y) a:= x \otimes ya$ defined for all $a, \ x, \ y \in A$ turn $A^{\otimes 2}$ into a $A$-bimodule. The $k$-vector space $\textsf{End}(A)$ of linear
endomorphisms of $A$ will be viewed as an associative algebra under composition denoted simply by concatenation $TS$, for $T,S \in \textsf{End}(A)$. The existence of a co-operation $\Delta: A \xrightarrow{} A^{\otimes 2}$ endows $\textsf{End}(A)$ with the so-called {\it{convolution product}} $*$, defined by $T *S := \mu(T \otimes S)\Delta$, for all $T,S \in \textsf{End}(A)$. It is associative if the co-operation $\Delta$ is coassociative.
\section{Algebraic structures related to one co-operation}
\subsection{$\epsilon'[R]$-Bialgebras and $L$-anti-dipterous algebras}
In this subsection, we show that associative structures can pop up from non-coassociative co-operations provided they well behave with the underlying associative product.
\begin{lemm}
\label{Lem AD}
Let $(A, \ \mu)$ be an associative algebra and $\zeta: A \xrightarrow{} A$ be a right Baxter operator. Define the binary operations, $\overleftarrow{\star}_\zeta, \ \prec_\zeta: A^{\otimes 2} \xrightarrow{} A$ by:
$$x \overleftarrow{\star}_\zeta y := \zeta(x)y, \ \ \ \ \ 
x \prec_\zeta y := x \zeta(y), \ \forall \ x,y \in A.$$
Then, $(A, \ \overleftarrow{\star}_\zeta, \ \prec_\zeta)$ is a $L$-anti-dipterous algebra.
\end{lemm}
\Proof
Let $(A, \ \mu)$ be an associative algebra and $\zeta: A \xrightarrow{} A$ be a right Baxter operator, i.e., for all $x,y \in A,$ 
$$ \zeta(x)\zeta(y) := \zeta(\zeta(x)y).$$
Let us show for instance that the binary operation $\overleftarrow{\star}_\zeta$ is associative. Indeed,
for all $x,y,z \in A$, $x \overleftarrow{\star}_\zeta (y \overleftarrow{\star}_\zeta z) := \zeta(x)\zeta(y)z$ and $(x \overleftarrow{\star}_\zeta y) \overleftarrow{\star}_\zeta z :=\zeta(\zeta(x)y)z$. Checking the two other axioms of Definition \ref{LAD} entails that $(A, \ \overleftarrow{\star}_\zeta, \ \prec_\zeta)$ is a $L$-anti-dipterous algebra.
\eproof

\noindent
One of the interest of $L$-anti-dipterous algebras can be found in the following theorem and its applications.
\begin{theo}
\label{ThAD}
Let $(A, \mu, \Delta)$ be a $\epsilon'[R]$-bialgebra. 
Equip the algebra $\textsf{End}(A)$ with the convolution product $*$, not supposed to be associative here.
Then, there exists a right Baxter operator, the right shift $\beta: \textsf{End}(A) \xrightarrow{} \textsf{End}(A)$ given by $\beta(T):= id * T$, for all $T \in \textsf{End}(A)$. Set,
$$T \overleftarrow{\star}_\beta S := \beta(T) S \ \ \ \ \textrm{and} \ \ \ \
T \prec_\beta S := T  \beta(S), \ \forall \ T,S \in \textsf{End}(A).$$
Then, $(\textsf{End}(A), \overleftarrow{\star}_\beta, \prec_\beta)$ is a $L$-anti-dipterous algebra.
\end{theo}
\Proof
Let $(A, \mu, \Delta)$ be a $\epsilon'[R]$-bialgebra. For all $x \in A$, write, using the Sweedler notation, $\Delta(x):= x_{(1)} \otimes x_{(2)}$.
Let $T, S \in \textsf{End}(A)$. On the one hand,
\begin{eqnarray*}
\beta(T)\beta(S)(x) &:=& \beta(T)(x_{(1)} S(x_{(2)})), \\
&=& \mu (id \otimes T)x_{(1)} S(x_{(2)})_{(1)} \otimes S(x_{(2)})_{(2)}, \\
&=& x_{(1)} S(x_{(2)})_{(1)} T( S(x_{(2)})_{(2)}).
\end{eqnarray*}
On the other hand,
\begin{eqnarray*}
\beta(\beta(T)S)(x) &:=& \mu(id \otimes \beta(T)S) (x_{(1)} \otimes x_{(2)}) = x_{(1)} \beta(T)(S(x_{(2)})), \\
&=& x_{(1)} S(x_{(2)})_{(1)} T( S(x_{(2)})_{(2)}),
\end{eqnarray*}
proving that the right shift $\beta$ is a right Baxter operator on $\textsf{End}(A)$.
Applying Lemma \ref{Lem AD} entails that $(\textsf{End}(A), \overleftarrow{\star}_\beta, \prec_\beta)$ is a $L$-anti-dipterous algebra.
\eproof
\Rk
From \textbf{non}-coassociative co-operations, associative algebras can be constructed via right Baxter operators. In terms of category, Theorem \ref{ThAD} yields,
$$\epsilon'[R]-\textsf{Bialg.} \longrightarrow \textsf{L-A.Dip.} \longrightarrow \textsf{As.}$$
\begin{prop}
\label{extendcop}
Let $As(S)$ be the free associative $k$-algebra generated by a non-empty set $S$. Fix a co-operation $\Delta: kS \xrightarrow{} kS^{\otimes 2}$ and extend it to a co-operation $\Delta_{\sharp}: As(S) \xrightarrow{} As(S)^{\otimes 2} $ such that for all words say $s_1 \ldots s_n$, $\Delta_{\sharp}(s_1 \ldots v_s):= s_1 \ldots s_{n-1} \Delta(s_n)$. 
Then, $(As(S),\Delta_{\sharp})$ is a $\epsilon'[R]$-bialgebra.
\end{prop}
\Proof
Straightforward by construction of the co-operation $\Delta_{\sharp}$.
\eproof
\begin{exam}{[Weighted directed graphs and dynamics]}
\label{graphdyn}
Let $G=(G_0, G_1, s, t)$ be a directed graph, supposed to be locally-finite, row-finite, without sink and source, equipped with a family of weights $(w_v)_{v \in G_0}$ and such that $s \times t: G_1 \xrightarrow{} G_0 \times G_0$ is injective. Consider the free $k$-vector space $kG_0$ spanned by $G_0$. Identify any directed arrow $v \longrightarrow w \in G_1$ with $v \otimes w$.
The set $G_1$ is then viewed as a subset of $kG_0^{\otimes 2}$. The family of weights $(w_v)_{v \in G_0}$ is then viewed as 
a family of maps $w_v: \ F_v \xrightarrow{} k$, where $F_v:=\{a \in G_{1}, \ s(a)=v \}$.
Define the co-operation $\Delta_M: kG_0 \xrightarrow{} kG_0^{\otimes 2}$ as follows \cite{Coa}:
$$\Delta_M (v) := \sum_{i: a_i \in F_v} \ w_v(a_i) \ v \otimes t(a_i),$$ for all $v \in G_0$. Extend it to $As(G_0)$ as in Proposition \ref{extendcop}. Then, $(As(G_0), \Delta_{M\sharp})$ is a $\epsilon'[R]$-bialgebra.
\end{exam}
\begin{exam}{[Substitutions and Language theory]}
Let $S$ be a non-empty set. For all $s \in S$, suppose there exit substitutions of the form $s \mapsto s^i 
_{1}s^i _{2}$, $i \in I_s$, which appear with probability $\mathbb{P}(s \mapsto s^i _{1}s^i _{2})$ and $\card(I_s) < \infty$ for all $s \in S$. Consider the free 
$k$-vector space $k S$ and define the co-operation 
$\Delta: kS  \xrightarrow{} k S^{\otimes 2}$ by,
$$\Delta (s) := \sum_{i \in I_s} \ \mathbb{P}(s \mapsto s^i _{1}s^i _{2}) \ s^i _{1} \otimes s^i _{2}.$$
Extend it to $As(S)$ as in Proposition \ref{extendcop}. The $k$-vector space $(As(S), \ \Delta_{\sharp})$ is then a $\epsilon'[R]$-bialgebra.

\noindent
All the possible dynamics of a string are then given by the operator $\mu \Delta$, where $\mu$ is the associative product representing the concatenation of two symbols.
For instance, consider the start symbol at time $t=0$, say $s$. Apply
the operator $\Delta$ and $\mu$ to obtain, 
$$\sum_{i \in I_s} \ \mathbb{P}(s \mapsto s^i _{1}s^i _{2}) \ s^i _{1} s^i _{2},$$
which is all the possible strings at time $t=1$. The probability to get the string, say $s^{i_0} _{1} s^{i_0} _{2}$, with $i_0 \in I_s$ at time $t=1$, is then 
$\mathbb{P}(s \mapsto s^{i_0} _{1}s^{i_0} _{2})$ and so forth. 
\end{exam}
The dynamics of a string has another algebraic structure since
from any $\epsilon'[R]$-bialgebras, another $L$-anti-dipterous algebra can be constructed.
\begin{prop}
\label{subas}
Let $(A, \ \mu, \ \Delta)$ be a $\epsilon'[R]$-bialgebra. 
Set $x \bowtie y := \mu(\Delta(x)) y$ and $x \prec_{A} y := x \mu(\Delta(y))$ defined for all $x,y \in A$.
Then, the $k$-vector space $(A, \ \bowtie, \ \prec_{A} )$ is a $L$-anti-dipterous algebra.
\end{prop}
\Proof
Straightforward by using $\mu(\Delta(xy))= x\mu(\Delta(y))$, for all $x,y \in A$.
\eproof

\noindent
Proposition \ref{subas} gives another functor, $\epsilon'[R]-\textsf{Bialg.} \longrightarrow \textsf{L-A.Dip.}$
\Rk
Denote by $AL$, the operad associated with $L$-anti-dipterous algebras. Then, the sequence $(AL(n))_{n \geq 1}$ starts with $1,2,5,14, \ldots$.
\subsection{On $\epsilon[R]$-bialgebras}
We suppose the co-operation $\Delta$ to be coassociative and study what kind of structures emerge.
\begin{lemm}
\label{Lem Dip}
Let $(A, \ \mu)$ be an associative algebra and $\theta: A \xrightarrow{} A$ be a left Baxter operator. Define the binary operations $\overrightarrow{\star}_\theta, \ \succ_\theta: A^{\otimes 2} \xrightarrow{} A,$ by:
$$x \overrightarrow{\star}_\theta y := x\theta(y), \ \ \ \ \ 
x \succ_\theta y := \theta(x) y, \ \forall \ x,y \in A.$$
Then, $(A, \ \overrightarrow{\star}_\theta, \ \succ_\theta)$ is a $L$-dipterous algebra.
\end{lemm}
\Proof
Straightforward, since by definition left Baxter operators $ \theta$ on $(A, \ \mu)$ obey the identity $\theta(x)\theta(y) := \theta(x \theta(y))$, for all $x,y \in A$.
\eproof

\noindent
As $\epsilon[R]$-bialgebras are also $\epsilon'[R]$-bialgebras, the operations $\overleftarrow{*_{\beta}}$ and $\prec_{\beta}$ defined in the previous subsection are related to $\overrightarrow{*_{\gamma}}$ and $\succ_{\gamma}$ via the following lemma.
\begin{lemm}
\label{Lem ench}
Let $(A, \ \mu)$ be an associative algebra and $\zeta, \ (\theta): A \xrightarrow{} A$ be right (left) Baxter operators.
Then,
$$ (x \overleftarrow{\star}_\zeta y) \overrightarrow{\star}_\theta z = x \overleftarrow{\star}_\zeta (y \overrightarrow{\star}_\theta z), $$
$$ (x \succ_\theta y) \prec_{\zeta} z = x \succ_\theta (y \prec_{\zeta} z), $$
for all $x,y \in A$.
\end{lemm}
\Proof
Straightforward.
\eproof
\begin{theo}
\label{tthdip}
Let $(A, \mu, \Delta)$ be a $\epsilon[R]$-bialgebra. Then, there exists a left Baxter operator, the left shift $\gamma: \textsf{End}(A) \xrightarrow{} \textsf{End}(A)$ given by $T \mapsto  T* id$. Set,
$$T \overrightarrow{\star}_\gamma S := T  \gamma(S), \ \ \textrm{and} \ \
T \succ_\gamma S := \gamma(T)  S, \ \forall \ T,S \in \textsf{End}(A).$$
Then, $(\textsf{End}(A), \overrightarrow{\star}_\gamma, \succ_\beta)$, is a $L$-dipterous algebra.

\noindent
Moreover, the operations described in Theorem \ref{ThAD} are related to the new operations by the following relations:
$$ (R \overleftarrow{\star}_\beta S) \overrightarrow{\star}_\gamma T = R \overleftarrow{\star}_\beta (S \overrightarrow{\star}_\gamma T),$$
$$ (R \succ_\gamma S) \prec_{\beta} T = R \succ_\gamma (S \prec_{\beta} T), $$
for all $R,S, T \in \textsf{End}(A)$.
\end{theo}
\Proof
Let $(A, \mu, \Delta)$ be a $\epsilon[R]$-bialgebra. Equip $\textsf{End}(A)$ with the convolution product $*$. Fix $T,S \in \textsf{End}(A)$, $x \in A$ and write $\Delta(x) := x_{(1)} \otimes x_{(2)}$. On the one hand,
\begin{eqnarray*}
\gamma(T)\gamma(S)(x) &:=& \gamma(T)(S(x_{(1)}) x_{(2)}), \\
& =& \mu (T \otimes id ) \ S(x_{(1)}) (x_{(2)})_{(1)} \otimes (x_{(2)})_{(2)}, \\
& =& T (S(x_{(1)}) (x_{(2)})_{(1)})  (x_{(2)})_{(2)} .
\end{eqnarray*}
One the other hand,
\begin{eqnarray*}
\gamma(T\gamma(S))(x) &:=& \mu (T\gamma(S) \otimes id) (x_{(1)} \otimes x_{(2)}), \\
&=& T\gamma(S) (x_{(1)} ) x_{(2)}, \\
&=& (T \mu(S \otimes id) (x_{(1)})_{(1)} \otimes (x_{(1)})_{(2)})x_{(2)}, \\
&=& T(S((x_{(1)})_{(1)})(x_{(1)})_{(2)})x_{(2)}.
\end{eqnarray*}
Therefore,
$$\gamma(T)\gamma(S) = \gamma(T\gamma(S)),$$
holds since the coproduct $\Delta$ is coassociative. The linear map $\gamma$ is thus a left Baxter operator. The other claims follow from Lemmas \ref{Lem Dip} and \ref{Lem ench}.
\eproof

\noindent
In terms of category, Theorem \ref{tthdip} gives, 
$$\epsilon[R]-\textsf{Bialg.}  \longrightarrow \textsf{L-Dip.} \longrightarrow \textsf{Dip.} \longrightarrow \textsf{As.}$$

\noindent
The coassociativity of the coproduct $\Delta$ will imply that the operators $\beta$ and $\gamma$ will commute. Let us show what kind of structures appear.
\begin{theo}
\label{RL Com}
Let $(A, \ \mu)$ be an associative algebra and $\zeta, \ (\theta): A \xrightarrow{} A$ be right (left) Baxter operators which commute, i.e., $\zeta \theta = \theta \zeta$. Then, the left Baxter operator $\theta$ is a left Baxter operator on the $L$-anti-dipterous algebra $(A, \ \overleftarrow{\star}_\zeta, \ \prec_\zeta)$, i.e., for all $x,y \in A$,
$$ \theta(x)\overleftarrow{\star}_\zeta \theta(y) = \theta(x \overleftarrow{\star}_\zeta \theta(y)), \ \ \ \theta(x) \prec_\zeta \theta(y) = \theta(x \prec_\zeta \theta(y)).$$
Moreover, this entails the existence of four extra binary operations $\bullet_1, \  \bullet_2,\  \bullet_3,\  \bullet_4: A^{\otimes 2} \xrightarrow{} A$  defined by (EQ. Th. \ref{RL Com}),
$$  x \bullet_1 y := x \overleftarrow{\star}_\zeta \theta(y) := \zeta(x)\theta(y), \ \ \ \ \   x \bullet_2 y := \theta(x) \overleftarrow{\star}_\zeta y := \theta\zeta(x)y, $$
 $$ x \bullet_3 y := \theta(x) \prec_\zeta y := \theta(x)\zeta(y),\ \ \ \ \
 x \bullet_4 y := x \prec_\zeta \theta(y) := x\theta\zeta(y),$$
verifying the matrix of relations $([M]_{ij})_{(i,j=1,\ldots 3)}$:
$$
(x \bullet_4 y)\bullet_4 z = x \bullet_4(y \bullet_1 z); \ \ \ \
(x \bullet_3 y)\bullet_4 z = x \bullet_3(y \bullet_1 z); \ \ \ \
(x \bullet_1 y)\bullet_4 z = x \bullet_1(y \bullet_4 z); \ \ \ \  
$$
$$
(x \bullet_4 y)\bullet_3 z= x \bullet_3(y \bullet_2 z); \ \ \ \
(x \bullet_2 y)\bullet_4 z= x \bullet_2(y \bullet_4 z); \ \ \ \
(x \bullet_1 y)\bullet_1 z= x \bullet_1(y \bullet_1 z); \ \ \ \
$$
$$
(x \bullet_1 y)\bullet_2 z = x \bullet_2(y \bullet_2 z); \ \ \ \
(x \bullet_1 y)\bullet_3 z = x \bullet_2(y \bullet_3 z); \ \ \ \ 
(x \bullet_2 y)\bullet_1 z = x \bullet_2(y \bullet_1 z), \ \ \ \ $$
for all $x,y,z \in A$. We shall say that $(A, \bullet_1, \  \bullet_2,\  \bullet_3,\  \bullet_4)$ has a $([M]_{ij})_{(i,j=1,\ldots 3)}$-algebra structure if its matrix of relations verifies $([M]_{ij})_{(i,j=1,\ldots 3)}$, up to a permutation of $S_9$.
\end{theo}
\Proof
Let $(A, \ \mu)$ be an associative algebra and $\zeta, \ (\theta): A \xrightarrow{} A$ be right (left) Baxter operators such that $\zeta \theta = \theta \zeta$. Fix $x,y,z \in A$. We have,
\begin{eqnarray*}
\theta(x)\overleftarrow{\star}_\zeta \theta(y) &=&\zeta( \theta(x)) \theta(y), \\
&=&\theta(\zeta(x)) \theta(y), \\
&=&\theta(\zeta(x) \theta(y)), \\
&=& \theta(x \overleftarrow{\star}_\zeta \theta(y)).
\end{eqnarray*}
Proceeding similarly for the operation $\prec_\zeta$ shows that $\theta$ is a left Baxter operator on $(A, \ \overleftarrow{\star}_\zeta, \ \prec_\zeta)$.

\noindent
Let us check for instance $[M]_{11}$, i.e., the anti-dipterous relation: $(x \bullet_4 y)\bullet_4 z = x \bullet_4(y \bullet_1 z)$. 
\begin{eqnarray*}
(x \bullet_4 y)\bullet_4 z &:=& (x \prec_\zeta \theta(y)) \prec_\zeta \theta(z), \\
&:=& x \prec_\zeta (\theta(y) \overleftarrow{\star}_\zeta \theta(z)), \\
&:=& x \prec_\zeta \theta(y \overleftarrow{\star}_\zeta \theta(z)), \\
&:=& x \bullet_4 (y \bullet_1 z).
\end{eqnarray*}
Proceeding similarly for the 8 other equations shows that $(A, \bullet_1, \  \bullet_2,\  \bullet_3,\  \bullet_4)$ has a $([M]_{ij})_{(i,j)}$-algebra structure.
\eproof

\noindent
Keep notation of Theorem \ref{RL Com} and observe that $(A, \bullet_1, \bullet_4)$ is a $L$-anti-dipterous algebra whereas $(A, \bullet_1, \bullet_2)$ is a $L$-dipterous algebra.
Therefore, in terms of category, Theorem \ref{RL Com} gives the following commutative diagram.
\begin{center}
$
\begin{array}{cccccc}
& &  &  \textsf{L-Dip.} & \\
& & \nearrow  & & \searrow  \\ 
\epsilon[R]-\textsf{Bialg.} & \longrightarrow &  ([\textsf{M}]_{(ij)}) & \longrightarrow & & \textsf{As.}  \\ 
& & \searrow &  &  \nearrow \\
& &  & \textsf{L-A.Dip.} &   \\ 
\end{array} 
$
\end{center}
\Rk
In Theorem \ref{RL Com} the r\^oles play by $\zeta$ and $\theta$ can be inversed, $\overleftarrow{\star}_\zeta$ becoming $\overrightarrow{\star}_\theta$ and $\prec_\zeta$ becomming $\succ_{\theta}$. Hence, $\zeta$ becomes a right Baxter operator on the $L$-dipterous algebra $(A,  \overrightarrow{*_{\theta}}, \succ_{\theta})$. In this case we will obtain new four binary operations,
$\tilde{\bullet}_1, \  \tilde{\bullet}_2,\  \tilde{\bullet}_3,\  \tilde{\bullet}_4: A^{\otimes 2} \xrightarrow{} A$  defined by,
$$  x \tilde{\bullet}_1 y := x \overrightarrow{\star}_\theta \zeta(y), \ \ \ \ \   x \tilde{\bullet}_2 y := \zeta(x) \overrightarrow{\star}_\theta y, $$
 $$ x \tilde{\bullet}_3 y := \zeta(x) \succ_\theta y, \ \ \ \ \
 x \tilde{\bullet}_4 y := x \succ_\theta \zeta(y).$$
However, observe that the new matrix of relations can be obtained from (EQ. Th. \ref{RL Com}) by the permutation $\sigma$ defined by $\sigma(1)=2, \ \sigma(2)=3, \ \sigma(3)=4, \ \sigma(4)=1$. For instance, $[M]_{11}$ becomes,
$$[(x \tilde{\bullet}_{\sigma(4)} y)\tilde{\bullet}_{\sigma(4)} z = x \tilde{\bullet}_{\sigma(4)}(y \tilde{\bullet}_{\sigma(1)} z)] := [(x \tilde{\bullet}_{1} y)\tilde{\bullet}_{1} z = x \tilde{\bullet}_{1}(y \tilde{\bullet}_{2} z)].$$ 
However, 
\begin{eqnarray*}
(x \tilde{\bullet}_{1} y)\tilde{\bullet}_{1} z &:=& (x \overrightarrow{\star}_\theta \xi(y))\overrightarrow{\star}_\theta \xi( z) = x \overrightarrow{\star}_\theta (\xi(y) \overrightarrow{\star}_\theta \xi(z)), \\
&=& x \overrightarrow{\star}_\theta \xi((\xi(y) \overrightarrow{\star}_\theta z)),\\
&=& x \tilde{\bullet}_{1}(y \tilde{\bullet}_{2} z).
\end{eqnarray*}
\begin{prop}{\bf{[Opposite structure of $([M]_{ij})_{(i,j)}$-algebras]}} 
\label{opo}
Let $(A, \bullet_1, \  \bullet_2,\  \bullet_3,\  \bullet_4)$ be a $([M]_{ij})_{(i,j)}$-algebra.
Define new binary operations $\bullet^{op}_1, \ \bullet^{op}_2, \ \bullet^{op}_3, \ \bullet^{op}_4: A^{\otimes 2} \xrightarrow{} A$ by,
$$x \bullet^{op}_1 y := y \bullet_1 x, \ \ \ x \bullet^{op}_3 y := y \bullet_3 x, \ \ \ x \bullet^{op}_2 y := y \bullet_4 x, \ \ \
x \bullet^{op}_4 y := y \bullet_2 x,$$ 
for all $x,y \in A$. Then,
$(A, \ \bullet^{op}_1, \ \bullet^{op}_2, \ \bullet^{op}_3, \ \bullet^{op}_4)$ is still a $([M]_{ij})_{(i,j=1,\ldots 3)}$-algebra.
\end{prop}
\Proof
Consider the matrix $( [M]_{ij})_{(i,j)}$ stated in Theorem \ref{RL Com}. Denote by $( [M]^{op}_{ij})_{(i,j)}$
 the matrix obtained from
$([M]_{ij})_{(i,j)}$
 by replacing $ [M]_{ij}$ by the transformation stated in Proposition \ref{opo}.
Observe the exitence of an axial symmetry around the row 2. Indeed,
for all $j=1, \ldots, 3$ $[M]_{2j} := [M]^{op}_{2j}$,
$[M]_{1j} := [M]^{op}_{3j}$ and $[M]_{3j} := [M]^{op}_{1j}$.
\eproof

\noindent
A $([M]_{ij})_{(i,j)}$-algebra
$(A, \ \bullet^{op}_1, \ \bullet^{op}_2, \ \bullet^{op}_3, \ \bullet^{op}_4)$ is said to be {\it{commutative}} when it coincides with its opposite structure, i.e., $x \bullet_1 y= y \bullet_1 x$, $x \bullet_2 y := y \bullet_4 x$ and $x \bullet_3 y= y \bullet_3 x$.
\Rk
The operad $P_{( [M]_{ij})_{(i,j)}}$
 associated with algebras of type $( [M]_{ij})_{(i,j)}$
is binary, quadratic and non-symmetric. The sequence $(\dim \ P_{( [M]_{ij})_{(i,j)}}(n))_{n \geq 1}$ starts with $1, \ 4, \ 23, \ 160 \ldots.$
\begin{coro}
Let $(A, \ \mu, \Delta)$ be a $\epsilon[R]$-bialgebra. Then, there exists a structure of $( [M]_{ij})_{(i,j)}$-algebra on $\textsf{End}(A)$.
\end{coro}
\Proof
Keep notation of Theorem \ref{RL Com} and 
apply it to the right Baxter operator $\beta$ which commutes with the left Baxter operator $\gamma$, where $\beta(T):= id * T$ and $\gamma(T):= T *id$, for all $T \in \textsf{End}(A)$.
\eproof
\begin{coro}
\label{Transp RR}
Let $(A, \ \mu)$ be an associative algebra and $\zeta, \ (\theta): A \xrightarrow{} A$ be two right Baxter operators which commute, i.e., $\zeta \theta = \theta \zeta$. Then, the right Baxter operator $\theta$ is a right Baxter operator on the $L$-anti-dipterous algebra $(A, \ \overleftarrow{\star}_\zeta, \ \prec_\zeta)$, i.e., for all $x,y \in A$,
$$ \theta(\theta(x)\overleftarrow{\star}_\zeta y) = \theta(\theta(x) \overleftarrow{\star}_\zeta y), \ \ \ \theta(x) \prec_\zeta \theta(y) = \theta(\theta(x) \prec_\zeta y).$$
Moreover, this entails the existence of four extra binary operations $\bullet_1, \  \bullet_2,\  \bullet_3,\  \bullet_4: A^{\otimes 2} \xrightarrow{} A$, defined in Theorem \ref{RL Com}, and verifying the matrix $([M_1]_{ij})_{(i,j=1,\ldots 3)}$,
$$
(x \bullet_2 y)\bullet_2 z = x \bullet_2(y \bullet_2 z); \ \ \ \
(x \bullet_1 y)\bullet_4 z = x \bullet_1(y \bullet_3 z); \ \ \ \ 
(x \bullet_2 y)\bullet_1 z = x \bullet_2(y \bullet_1 z); \ \ \ \  
$$
$$
(x \bullet_3 y)\bullet_4 z = x \bullet_3(y \bullet_1 z); \ \ \ \
(x \bullet_2 y)\bullet_4 z= x \bullet_2(y \bullet_4 z); \ \ \ \
(x \bullet_1 y)\bullet_1 z= x \bullet_1(y \bullet_2 z); \ \ \ \
$$
$$
(x \bullet_2 y)\bullet_3 z = x \bullet_2(y \bullet_3 z); \ \ \ \
(x \bullet_3 y)\bullet_3 z= x \bullet_3(y \bullet_2 z); \ \ \ \
(x \bullet_4 y)\bullet_4 z =x  \bullet_4(y \bullet_2 z). \ \ \ \
$$
We shall say that the $k$-vector space $(A, \ \bullet_1, \  \bullet_2,\  \bullet_3,\  \bullet_4 )$ is an algebra of type $([M_1]_{ij})_{(i,j=1,\ldots 3)}$ if its matrix of relations is $([M_1]_{ij})_{(i,j=1,\ldots 3)}$ up to a permutation of $S_9$.
\end{coro}
\Proof
As the Baxter operator $\theta$ is no longer left but right, this implies the existence of the following relations,
$$ x \bullet_1(y \bullet_i z) \longleftrightarrow x \bullet_1(y \bullet_j z), \  x \bullet_4(y \bullet_i z) \longleftrightarrow x \bullet_4(y \bullet_j z), $$ $$  (x \bullet_i y)\bullet_3 z \longleftrightarrow (x \bullet_j y)\bullet_3 z,  \ (x \bullet_i y)\bullet_2 z \longleftrightarrow (x \bullet_j y)\bullet_2 z,$$
with if $i=1$, (resp. $i=4$) then $j$ has to be removed by $2$ (resp. by 3). Apply these relations to the matrix $( [M]_{ij})_{(i,j=1,\ldots 3)}$
of Theorem \ref{RL Com} to obtain the matrix $([M_1]_{ij})_{(i,j=1,\ldots 3)}$ up to a permutation of $S_9$. 
\eproof

\noindent
Keep notation of Corollary \ref{Transp RR} and observe that both $(A, \bullet_2, \bullet_3)$ and $(A, \bullet_2, \bullet_1)$ are $L$-anti-dipterous algebras.
Therefore, Corollary \ref{Transp RR} gives the following commutative diagram.
\begin{center}
$
\begin{array}{cccccc}
& &  &  \textsf{L-A.Dip.} & \\
& & \nearrow  & & \searrow  \\ 
 & &  ([\textsf{M}_1]_{(ij)}) & \longrightarrow & & \textsf{As.}  \\ 
& & \searrow &  &  \nearrow \\
& &  & \textsf{L-A.Dip.} &   \\ 
\end{array} 
$
\end{center}
\Rk
Denote by $P_{([M_1]_{ij})_{(i,j)}}$, the binary quadratic non-symmetric operad associated with algebra of type $([M_1]_{ij})_{ (i,j=1, \ldots 3) }$.
The sequence $(P_{([M_1]_{ij})_{(i,j)}}(n))_{n \geq 1}$ starts
with $1,4,23,160 \ldots$.
\begin{prop} {\bf{[Transpose of algebras of type $([M_1]_{ij})_{(i,j)}$
]}} 
Let $(A, \ \bullet_{1}, \ \bullet_{2},\ \bullet_{3}, \ \bullet_{4}) $ be an algebra of type $([M_1]_{ij})_{(i,j)}$. Define for all $x,y \in A$ four binary operations, $\bullet_{1,t}, \ \bullet_{2,t},\ \bullet_{3,t}, \ \bullet_{4,t}: A^{\otimes 2} \xrightarrow{} A $ by,
$$ (Tr \ 2.) \ \ x \bullet_{1,t}  y := x \bullet_{3}  y, \ \
 x \bullet_{3,t}  y := x \bullet_{1}  y, \ \
x \bullet_{2,t}  y := x \bullet_{2}  y, \ \
x \bullet_{4,t}  y := x \bullet_{4}  y.$$
Denote by $([M_1]^t_{ij})_{(i,j)}$
the matrix of relations obtained from $([M_1]_{ij})_{(i,j)}$ by applying the transformation $(Tr 2.)$.
Then, the $k$-vector space $A$ equipped with these four binary operations is still an algebra of type $([M_1]_{ij})_{(i,j)}$.
\end{prop}
\Proof
Observe that the transformation $(Tr \ 2.)$ let the matrix $([M_1]_{ij})_{(i,j=1,\ldots 3)}$
 globally invariant since for all $i,j=1,\ldots 3$,
$[M_1]_{ij} = [M_1]_{ji}.$
\eproof
\Rk
Another extra algebraic structure can be exhibited.
Let $(A, \mu, \Delta)$ be a $\epsilon[R]$-bialgebra.
There exist 3 extras binary operations $\searrow,\nearrow,
\nwarrow: A^{\otimes 2} \xrightarrow{} A$, given by, 
$$x \nearrow y := x \overleftarrow{\star}_\beta \ \gamma(y):= \beta(x) \overrightarrow{\star}_\gamma y, \ \ \ \ \ \ x \nwarrow y:= x \prec_{\beta}  \ \gamma(y), \ \ \ \ \ \ x \searrow y:= \beta(x) \succ_{\gamma} y,$$
and verifying:
\begin{eqnarray*}
& &   (x \nearrow y)\nearrow z = x \nearrow(y \nearrow z), \ \ \
(x \nwarrow y)\nwarrow z = x \nwarrow(y \nearrow z),\ \ \
  (x \searrow y)\nearrow z =  x \searrow(y \nearrow z), \\
& & [  (x \nearrow y)\nearrow z = x \nearrow(y \nearrow z)], \ \ \ (x \nearrow y)\searrow z = x \searrow(y \searrow z),\ \ \ (x \nearrow y)\nwarrow z = x \nearrow(y \nwarrow z), \\
& & (x \searrow y)\nwarrow z = x \searrow(y \nwarrow z).
\end{eqnarray*}
\Rk
Observe that the first row describes $L$-anti-dipterous algebra
axioms, the second one $L$-dipterous algebra axioms, these two structures being linked by the last relation. The sequence of dimensions of the associated operad starts with $1, 3, 12, \ldots$.

\noindent
Let us construct examples of $\epsilon[R]$-bialgebras.
\begin{exam}{}
Fix $\lambda \in k$.
Let $A$ be the unital commutative $k$-algebra $k[X_1, \ldots, X_n] /  \mathcal{R}$, where $\mathcal{R}$ is the set of relations $X_iX_j =0$, for all $i = 1, \ldots n$. Set $\Sigma:= \sum_{i=1, \ldots,n} \ X_i$. Then, the following co-operation $\Delta_{\lambda}:  A \xrightarrow{} A^{\otimes 2}$, defined by $\Delta_{\lambda}(1):= 1 \otimes \Sigma + \Sigma \otimes 1 + \lambda \Sigma \otimes \Sigma$ and $\Delta_{\lambda}(X_i):= X_i \otimes \Sigma$ is coassociative and turns $A$ into a $\epsilon[R]$-bialgebra. There exists a linear map $\epsilon_{\lambda}: A \xrightarrow{} k$ given by
$\epsilon_{\lambda}(\Sigma):= 1$ and $\epsilon_{\lambda} (1):= -\lambda$ verifying $(id \otimes \epsilon_{\lambda})\Delta_{\lambda}(1) =1 =(\epsilon_{\lambda} \otimes id )\Delta_{\lambda}(1)$ and $(id \otimes \epsilon_{\lambda})\Delta_{\lambda}(X_i) =X_i$. Observe that in the unital commutative $k$-algebra $k[X]$, the coproduct is cocommutative. Such a structure appears in \cite{Khovanov}.
\end{exam}
\begin{prop}
Let $S$ be a non-empty set and $\Delta: kS \xrightarrow{} kS^{\otimes 2}$ be a coassociative co-operation. Consider the free associative $k$-algebra $As(S)$ generated by $S$. Extend the coproduct $\Delta$ to $As(S)$ by considering $\Delta_{\sharp}: As(S) \xrightarrow{} As(S)^{\otimes 2}$ defined by $\Delta_{\sharp}(s_1 \ldots s_n):= s_1 \ldots s_{n-1} \Delta(s_n)$, where $s_1, \ldots, s_n \in S$. Then, $(As(S), \ \Delta_{\sharp})$ is a $\epsilon[R]$-bialgebra.
If $(\eta \otimes id)\Delta=id$ (resp. $(id \otimes \eta)\Delta=id$) holds for a linear map
$\eta: kS \xrightarrow{} k$ then so are the equations $(\eta_{\sharp} \otimes id)\Delta_{\sharp}=id$ (resp. $(id \otimes \eta_{\sharp})\Delta_{\sharp}=id$)
for the map $\eta_{\sharp}: As(S) \xrightarrow{} As(S) \oplus k$ given by $\eta_{\sharp}(s_1 \ldots s_{n-1}s_n) := s_1 \ldots s_{n-1} \eta(s_n)$, where $s_1, \ldots, s_n \in S$. 
\end{prop}
\Proof
Let $S$ be a set and $\Delta: kS \xrightarrow{} kS^{\otimes 2}$ be a coassociative co-operation. Consider the free associative $k$-algebra $As(S)$ generated by $S$. We have to show that $\Delta_{\sharp}: As(S) \xrightarrow{} As(S)^{\otimes 2}$ is coassociative. Fix $s_1, \ldots s_n \in S$. As $s_1 \ldots s_n \in As(S)$ and $\Delta(s_n) := (s_n)_{(1)} \otimes (s_n)_{(2)}$ we get:
$$(id \otimes \Delta_{\sharp}) \Delta_{\sharp}(s_1 \ldots s_n):=s_1 \ldots s_{n-1}(s_n)_{(1)} \otimes \Delta((s_n)_{(2)}),$$
and 
$$(\Delta_{\sharp} \otimes id) \Delta_{\sharp}(s_1 \ldots s_n):=s_1 \ldots s_{n-1}\Delta((s_n)_{(1)}) \otimes (s_n)_{(2)}.$$
Therefore, $(As(S), \ \Delta_{\sharp})$ is a $\epsilon[R]$-bialgebra.
\eproof
\Rk
As an application, we can deal with language theory and coassociative substitutions, see \cite{perorb1, dipt}. Let $S$ be a non-empty set.  Suppose substitutions on $S$ are of the form $s \mapsto s_{1}^i s_{2}^i$, with $i \in I_s$ and $\card I_s < \infty$. Consider the free $k$-vector space $kS$ spanned by $S$ and construct the co-operation $\Delta: s \mapsto \sum_{i \in I_s}  \ s_{1}^i \otimes s_{2}^i$. If $\Delta$ is coassociative, then the set of substitutions on $S$ is said to be {\it{coassociative}}. Such substitutions appear in an algebraic description of the chaotic map $x \mapsto 2x \mod 1$, with $x \in [0,1]$ and a quantisation of the Bernoulli walk \cite{perorb1}. They also appear in works on coassociative tilings and coverings of directed graphs \cite{dipt,codialg1}.
\begin{exam}{[$M_2(k)$]} 
\label{$M_2(k)$} 
We give now another method, inspired from \cite{Aguiar} to obtain $\epsilon[R]$-bialgebras.
\begin{prop}
Let $(A, \ \mu)$ be an associative algebra and $J$ be a finite set. Fix $A_j,B_j \in A$, $j \in J$. Consider the co-operation $\Delta: A \xrightarrow{} A^{\otimes 2}$ defined by $x \mapsto \Delta(x):= \sum_j x A_j \otimes B_j$.
Then, $(A, \ \mu, \ \Delta)$ is a $\epsilon[R]$-bialgebra if,
$$ \sum_{i,j} A_i \otimes B_iA_j \otimes B_j = \sum_{i,j} A_iA_j\otimes B_j \otimes B_i.$$ 
\end{prop}
\Proof
Straightforward.
\eproof

\noindent
Let $(A, \ \mu)$ be an associative algebra and $\lambda \in k$. Consider the co-operation $\Delta_1: A \xrightarrow{} A^{\otimes 2}$ defined by $x \mapsto \Delta_1(x):=  x A_1 \otimes B_1$. Then the co-operation $\Delta_1$ is coassociative if $A_1 ^2= \lambda A_1$ and $B_1A_1=\lambda B_1$. For instance consider in $M_2(k)$ the following matrices:
\[A = \begin{pmatrix}
 \lambda & \lambda\\
0 & 0
\end{pmatrix}, \ \ \ \
B = \begin{pmatrix}
 0 & 0 \\
\lambda & \lambda
\end{pmatrix}.
\]
Equipped with this coproduct, $M_2(k)$ is a  $\epsilon[R]$-bialgebra.
\end{exam}
\subsection{$\epsilon'[L]$-Bialgebras and $L$-anti-dipterous algebras}
What was done with $\epsilon[R]$-bialgebras can be realised with their left versions. The r\^ole played by the right Baxter operator $\beta$ is now hold in the left version by the operator $\gamma$ and conversely. 
\begin{theo}
\label{Leftadip}
Let $(A, \mu, \Delta)$ be a $\epsilon'[L]$-bialgebra. 
Equip the $k$-algebra $\textsf{End}(A)$ with the convolution product $*$, not supposed to be associative here.
Then, the left shift $\gamma: \textsf{End}(A) \xrightarrow{} \textsf{End}(A)$ given by $T \mapsto T * id$ is a right Baxter operator. Set,
$$T \overleftarrow{\star}_\gamma S := \gamma(T)  S \ \ \ \ \textrm{and} \ \ \ \
T \prec_\gamma S := T  \gamma(S), \ \ \forall T,S \in \textsf{End}(A).$$
Then, $(\textsf{End}(A), \overleftarrow{\star}_\gamma, \prec_\gamma)$ is a $L$-anti-dipterous algebra.
\end{theo}
\Proof
Let $(A, \mu, \Delta)$ be a $\epsilon'[L]$-bialgebra. Fix $T, S \in \textsf{End}(A)$, $x \in A$ and set $\Delta(x):= x_{(1)} \otimes x_{(2)}$. On the one hand,
\begin{eqnarray*}
\gamma(T)\gamma(S)(x) &:=& \gamma(T)(S(x_{(1)}) x_{(2)}), \\
&=& \mu (T \otimes id) S(x_{(1)})_{(1)} \otimes S(x_{(1)})_{(2)}x_{(2)}, \\
&=& T(S(x_{(1)})_{(1)}) (S(x_{(1)})_{(2)}x_{(2)}.
\end{eqnarray*}
On the other hand,
\begin{eqnarray*}
\gamma(\gamma(T)S)(x) &:=& \mu(\gamma(T)S \otimes id) (x_{(1)} \otimes x_{(2)}) = \gamma(T)S(x_{(1)}) x_{(2)}, \\
&=&  T(S(x_{(1)})_{(1)}) (S(x_{(1)})_{(2)} x_{(2)},
\end{eqnarray*}
proving that $\gamma$ is a right Baxter operator.
Applying Lemma \ref{Lem AD} entails that $(\textsf{End}(A), \overleftarrow{\star}_\gamma, \prec_\gamma)$ is a $L$-anti-dipterous algebra.
\eproof

\noindent
In terms of category, we get: 
$$\epsilon'[L]-\textsf{Bialg.} \longrightarrow \textsf{L-A.Dip.} \longrightarrow \textsf{As.}$$
Let us construct examples of $\epsilon'[L]$-bialgebras. 
\begin{prop}
\label{Pr LAN}
Let $As(S)$ be the free associative $k$-algebra generated by a non-empty set $S$. Fix a co-operation $\Delta: kS \xrightarrow{} kS^{\otimes 2}$ and extend it to a co-operation $\Delta_{\sharp}: As(S) \xrightarrow{} As(S)^{\otimes 2} $ such that for all words say $s_1 \ldots s_n$, $\Delta_{\sharp}(s_1 \ldots s_n):= \Delta(s_1)s_2 \ldots s_{n-1} s_n$. 
Then, $(As(S),\Delta_{\sharp})$ is a $\epsilon'[L]$-bialgebra.
\end{prop}
\Proof
Straightforward by construction of the co-operation $\Delta_{\sharp}$.
\eproof
\begin{exam}{[Weighted directed graphs and dynamics]}
Keep notation of example \ref{graphdyn}.
Define the co-operation $\tilde{\Delta}_M : kG_0 \xrightarrow{} kG_0^{\otimes 2}$ as follows \cite{Coa}:
$$\tilde{\Delta}_M (v) := \sum_{i: a_i \in P_v} \ w_{s(a_i)}(a_i) \ s(a_i) \otimes v,$$
 where $P_{v} := \{a \in G_{1}, \ t(a)=v \}$, for all $v \in G_0$. Then, $(As(G_0), \tilde{\Delta}_{\sharp M})$ is a $\epsilon'[L]$-bialgebra.
\end{exam}
More generally, we get,
\begin{exam}{[Substitutions and Language theory]}
Let $S$ be a non-empty set. For all $s \in S$, suppose there exit substitutions of the form $s \mapsto s^i _{1}s^i _{2}$, $i \in I_s$ with $\card(I_s) < \infty$, which appear with probablity $\mathbb{P}(s \mapsto s^i _{1}s^i _{2})$. Consider the free 
$k$-vector space $k S$ spanned by $S$ and define the co-operation 
$\Delta: kS  \xrightarrow{} k S^{\otimes 2}$ by,
$$\Delta (s) := \sum_{i \in I_s} \ \mathbb{P}(s \mapsto s^i _{1}s^i _{2}) \ s^i _{1} \otimes s^i _{2}.$$
Then, extend it to $As(S)$ as in Proposition \ref{Pr LAN}. The $k$-vector space $(As(S), \ \Delta_{\sharp})$ is a $\epsilon'[L]$-bialgebra.
\end{exam}

\noindent
From any $\epsilon'[L]$-bialgebras, a $L$-dipterous algebra can be constructed. 
\begin{prop}
Let $(A, \ \mu, \ \Delta)$ be a $\epsilon'[L]$-bialgebra. 
Set $x \succ_{A} y := \mu(\Delta(x)) y$ and $x \bowtie y := x \mu(\Delta(y))$ defined for all $x,y \in A$.
Then, the $k$-vector space $(A, \ \bowtie, \ \succ_{A} )$ is a $L$-dipterous algebra.
\end{prop}
\Proof
Straightforward by using $\mu(\Delta(xy))= \mu(\Delta(x))y $, for all $x,y \in A$.
\eproof
\Rk
Applied this proposition to model substitutions in symbolic dynamics via the associative product $\bowtie$ or the operation $\succ_{A}$.
\subsection{On $\epsilon[L]$-bialgebras}
In this subsection, the co-operation $\Delta$ is supposed to be coassociative.
\begin{theo}
Let $(A, \mu, \Delta)$ be a $\epsilon[L]$-bialgebra. Then, the right shift $\beta: \textsf{End}(A) \xrightarrow{} \textsf{End}(A)$ given by $T \mapsto  id *T$ is a left Baxter operator. Set,
$$T \overrightarrow{\star}_\beta S := T \beta(S), \ \ \textrm{and} \ \
T \succ_\beta S := \beta(T) S, \ \ \forall T, S \in \textsf{End}(A).$$
Then, $(\textsf{End}(A), \overrightarrow{\star}_\beta, \succ_\beta)$ is a $L$-dipterous algebra. Therefore there exists a functor such that,
$$\epsilon[L]-\textsf{Bialg.}  \longrightarrow \textsf{L. Dip.}$$

\noindent
These two operations are related to $\prec_{\gamma}$ and $\overleftarrow{\star}_\gamma$ from Theorem \ref{Leftadip} by,
$$ (R \overleftarrow{\star}_\gamma S) \overrightarrow{\star}_\beta T = R \overleftarrow{\star}_\gamma (S \overrightarrow{\star}_\beta T),$$
$$ (R \succ_\beta S) \prec_{\gamma} T = R \succ_\beta (S \prec_{\gamma} T), $$
for all $R,S, T \in \textsf{End}(A)$.
\end{theo}
\Proof
Let $(A, \mu, \Delta)$ be a $\epsilon[L]$-bialgebra. Equip $\textsf{End}(A)$ with the convolution product $*$. Fix $T,S \in \textsf{End}(A)$, $x \in A$ and set $\Delta(x) := x_{(1)} \otimes x_{(2)}$. On the one hand,
\begin{eqnarray*}
\beta(T)\beta(S)(x) &:=& \beta(T)(x_{(1)} S(x_{(2)})), \\
& =& \mu (id \otimes T ) \ (x_{(1)})_{(1)} \otimes (x_{(1)})_{(2)}S(x_{(2)}), \\
& =& (x_{(1)})_{(1)} T( (x_{(1)})_{(2)}S(x_{(2)})).
\end{eqnarray*}
One the other hand,
\begin{eqnarray*}
\beta(T \beta(S))(x) &:=& \mu (id \otimes T\beta(S)) (x_{(1)} \otimes x_{(2)}), \\
&=& x_{(1)}  T\beta(S)( x_{(2)}), \\
&=& x_{(1)} T \mu(id \otimes S) (( x_{(2)})_{(1)} \otimes (x_{(2)})_{(2)}),\\
&=& x_{(1)} T (( x_{(2)})_{(1)} S( (x_{(2)})_{(2)})).
\end{eqnarray*}
Since the coproduct $\Delta$ is coassociative, 
$$\beta(T)\beta(S) = \beta(T\beta(S)),$$ holds. The linear map $\beta$ is a left Baxter operator. The other claims follow from Lemmas \ref{Lem Dip} and \ref{Lem ench}.
\eproof
\Rk
In the left version, the r\^oles of the operators $\beta$ and $\gamma$ are reversed. 
\begin{coro}
Let $(A, \ \mu, \Delta)$ be a $\epsilon[L]$-bialgebra. Then, there exists a $( [M]_{ij})_{(i,j)}$-algebra structure on $\textsf{End}(A)$.
\end{coro}
\Proof
Keep notation of Theorem \ref{RL Com} and 
apply it to the right Baxter operator $\gamma$ which commutes with the left Baxter operator $\beta$, where $\beta(T):= id * T$ and $\gamma(T):= T *id$, for all $T \in \textsf{End}(A)$.
\eproof

\noindent
In terms of category, the following diagram commutes.
\begin{center}
$
\begin{array}{cccccc}
& &  &  \textsf{L-Dip.} & \\
& & \nearrow  & & \searrow  \\ 
\epsilon[L]-\textsf{Bialg.} & \longrightarrow &  ([\textsf{M}]_{ij}) & \longrightarrow & & \textsf{As.}  \\ 
& & \searrow &  &  \nearrow \\
& &  & \textsf{L-A.Dip.} &   \\ 
\end{array} 
$
\end{center}
\begin{prop}
\label{Prop Dip}
Let $(A, \ \mu)$ be an associative algebra and $\zeta, \ \theta: A \xrightarrow{} A$ be two left Baxter operators which commute, i.e., $\zeta \theta = \theta \zeta$. Then, the left Baxter operator $\theta$ is a left Baxter operator on the $L$-dipterous algebra $(A, \ \overrightarrow{\star}_\zeta, \ \succ_\zeta)$, i.e., for all $x,y \in A$,
$$ \theta(x)\overrightarrow{\star}_\zeta \theta(y) = \theta(x \overrightarrow{\star}_\zeta \theta(y)), \ \ \ \theta(x) \succ_\zeta \theta(y) = \theta(x \succ_\zeta \theta(y)).$$
Moreover, this entails the existence of four extra binary operations $\bullet_1, \  \bullet_2,\  \bullet_3,\  \bullet_4: A^{\otimes 2} \xrightarrow{} A$, defined by,
$$  x \bullet_4 y := x \overrightarrow{\star}_\zeta \theta(y) := x\zeta(\theta(y)), \ \ \ \ \   x \bullet_3 y := \theta(x) \overrightarrow{\star}_\zeta y := \theta(x)\zeta(y), $$
$$ x \bullet_2 y := \theta(x) \succ_\zeta y := \zeta(\theta(x))y,\ \ \ \ \
 x \bullet_1 y := x \succ_\zeta \theta(y) := \zeta(x)\theta(y),$$
verifying the matrix of relations $([M_2]_{ij})_{(i,j=1,\ldots 3)}$:
$$
(x \bullet_4 y)\bullet_4 z = x \bullet_4(y \bullet_4 z); \ \ \ \
(x \bullet_3 y)\bullet_1 z= x \bullet_2(y \bullet_1 z); \ \ \ \ 
(x \bullet_1 y)\bullet_4 z= x \bullet_1(y \bullet_4 z); \ \ \ \ $$
$$
(x \bullet_1 y)\bullet_3 z = x \bullet_2(y \bullet_3 z); \ \ \ \
(x \bullet_2 y)\bullet_4 z= x \bullet_2(y \bullet_4 z); \ \ \ \
(x \bullet_4 y)\bullet_1 z =x \bullet_1(y \bullet_1 z); \ \ \ \ $$
$$
(x \bullet_3 y)\bullet_4 z = x \bullet_3(y \bullet_4 z); \ \ \ \
(x \bullet_4 y )\bullet_3 z=x \bullet_3(y \bullet_3 z); \ \ \ \
(x \bullet_4 y)\bullet_2 z= x \bullet_2(y \bullet_2 z). \ \ \ \
$$
We shall say that the $k$-vector space $(A, \ \bullet_1, \  \bullet_2,\  \bullet_3,\  \bullet_4)$ is an algebra of type $([M_2]_{ij})_{(i,j)}$ if its matrix of relations is equal 
to $([M_2]_{ij})_{(i,j)}$ up to a permutation of $S_9$.
\end{prop}
\Proof
Straightforward.
\eproof

\noindent
Observe that both the $k$-vector spaces $(A, \ \bullet_4, \  \bullet_1)$ and $(A, \  \bullet_4, \  \bullet_3)$ are $L$-dipterous algebras. Therefore the following diagram,

\begin{center}
$
\begin{array}{cccccc}
& &  &  \textsf{L-Dip.} & \\
& & \nearrow  & & \searrow  \\ 
\\
 &  &  ([\textsf{M}_2]_{ij}) & \longrightarrow & & \textsf{As.}  \\ 
\\
& & \searrow &  &  \nearrow \\
& &  & \textsf{L-Dip.} &   \\ 
\end{array} 
$
\end{center}
commutes.
\begin{prop}{} [\textbf{Transpose of an algebra of type $([M_2]_{ij} )_{(i,j)}$}]
Keep notation of Proposition \ref{Prop Dip}. Let $A$ be a $([M_2]_{ij} )_{(i,j)}$-algebra. Define four new binary operations $\bullet_{i,t}$, $i =1, \ldots 4$ by,
$$ \bullet_{1,t} \mapsto \bullet_{3}, \ \ \ \bullet_{2,t} \mapsto \bullet_{2}, \ \ \ \bullet_{3,t} \mapsto  \bullet_{1}, \ \ \ \bullet_{4,t} \mapsto \bullet_{4}. $$
Then, $(A, \ \bullet_{1,t}, \ \bullet_{2,t}, \ \bullet_{3,t}, \ \bullet_{4,t})$ is still an algebra of type $ ([M_2]_{ij} )_{(i,j)} $.
\end{prop}
\Proof 
Denote by $([M_2]^t _{ij})_{i,j}$, the matrix of relations obtained from $([M_2]_{ij})_{i,j}$ by applying the said transformation. Then, observe that 
$[M_2]_{ij} := [M_2]^t_{ji}$, for all $i,j=1, \ldots 3$.
\eproof
\Rk
Keep notation of Proposition \ref{Prop Dip}. The r\^ole of the operators
$\gamma$ and $\xi$ can be exchanged. In this case, the $ ([M_2]_{ij} )_{(i,j)} $-algebra we obtain is the transpose of the first established in 
Proposition \ref{Prop Dip}.
\begin{prop}
\label{opll}
Let $A$ be a $([M_2]_{ij})_{(i,j)}$-algebra. Consider the following binary operations
$$ x \bullet_{1,op} y := y \bullet_{1} x, \ \ \ x \bullet_{2,op} y := y \bullet_{4} x, \ \ \ x \bullet_{3,op} y := y \bullet_{3} x, \ \ \ x \bullet_{4,op} y  :=  y \bullet_{2} x, $$
for all $x,y \in A$.
Then, $(A, \ \bullet_{1,op}, \  \bullet_{2,op}, \ \bullet_{3,op}, \ \bullet_{4,op})$ is a $([M_1]_{i,j})_{(i,j)}$-algebra.
\end{prop}
\Proof
Denote by $([M_2]^{op}_{ij})_{(i,j)}$, the matrix of relations obtained from $([M_2]_{ij})_{(i,j)}$ by applying the said transformation. Observe that:
$([M_2]^{op}_{ij} = [M_1]_{ij}$, for all $i,j=1, \ldots 3$, (and up to a permutation of $S_9$ if the matrices are written in a different order).
\eproof

\noindent
Proposition \ref{opll} can be re-written in terms of category. By oberving that the transformation $op$ is one-to-one, we get,
$$([\textsf{M}_2]_{ij})  \leftrightarrows ([\textsf{M}_1]_{ij}).$$

\subsection{$\epsilon$-Bialgebras}
An {\it{infinitesimal bialgebra}} (abbreviated $\epsilon$-bialgebra) is a triple $(A, \ \mu, \ \Delta )$ where $(A, \ \mu)$ is an associative algebra and $(A,\ \Delta )$ is a coassociative coalgebra such that
for all $a,b \in A$,
$\Delta (ab)= a_{(1)} \otimes a_{(2)} b + a b_{(1)} \otimes b_{(2)}$.
Such bialgebras
appeared for the first time in the work of Joni and Rota \cite{Rota}, see also \cite{Aguiar, AguiarLoday, Lertribax, Voi}.
$\epsilon$-Bialgebras allow the construction of quadri-algebras \cite{AguiarLoday}. For completeness, some definitions and results are recalled.

\begin{defi}{[Quadri-algebra]} 
\label{quadri}
A {\it{quadri-algebra}} $Q$ is a $k$-vector space equipped with four binary operations $\nearrow, \ \searrow, \ \swarrow, \ \nwarrow: A^{\otimes 2} \xrightarrow{} A$ obeying 9 relations. To state them, we introduce the following sums.
\begin{eqnarray*}
x \succ y &:=& x \nearrow y + x \searrow y, \\
x  \prec y &:=& x \nwarrow y + x \swarrow y, \\
x \vee y &:=& x \searrow y + x \swarrow y, \\
x \wedge y &:=& x \nearrow y + x \nwarrow y, \\
\end{eqnarray*}
and $$ x \bar{\star} y := x \prec y + x  \succ y = x \wedge y + x \vee y,$$
for all $x,y \in Q$.
The
compatibility axioms are presented in the following $3 \times 3$ matrix $([Q]_{ij})_{(i,j= 1, \ldots 3)}$. 

$
\begin{array}{ccccccccc}
(x \nwarrow y)\nwarrow z  &=& x \nwarrow (y \bar{\star} z); &
(x \nearrow y)\nwarrow z &=& x \nearrow (y \prec z); &
(x \wedge y)\nearrow z &=& x \nearrow (y \succ z); \\
(x \swarrow y)\nwarrow z &=& x \swarrow (y \wedge z); &
(x \searrow y)\nwarrow z &=& x \searrow (y \nwarrow z); &
(x \vee y)\nearrow z &=& x \searrow (y \nearrow z); \\
(x \prec y)\swarrow z &=& x \swarrow (y \vee z); &
(x \succ y)\swarrow z &=& x \searrow (y \swarrow z); &
(x \bar{\star} y)\searrow z &=& x \searrow (y \searrow z). \\
\end{array} 
$
\end{defi}
\begin{theo} [\cite{AguiarLoday}]
Let $(A, \ \mu, \ \Delta)$ be a $\epsilon$-bialgebra. Then, there exist two commuting Baxter operators $\beta$ and $\gamma$ and thus a quadri-algebra structure on $\textsf{End}(A)$ given by:
$$ T \swarrow S := \gamma(T) \beta(S), \ \ \
T \nwarrow S := T \gamma \beta(S), \ \ \ 
T \searrow S := \gamma \beta(T) S, \ \ \
T \nearrow S := \beta(T)\gamma(S),
$$
with $\beta(T):= id * T$ and $\gamma(T):= T * id$, for all $T, S \in \textsf{End}(A)$.
\end{theo}
We recall the existence of a notion of commutative quadri-algebras and of transpose of quadri-algebras. The free quadri-algebra could be related to non-crossing connected graphs \cite{AguiarLoday}. The augmented free quadri-algebra has a connected Hopf algebra structure \cite{Lodayscd}. All these terms are explained in Section \ref{Octo}. 
A quadri-algebra gives two dendriform dialgebras, the so-called vertical structure $Q_v :=(Q, \ \wedge, \ \vee)$ and the so-called horizontal structure $Q_h :=(Q, \ \prec, \ \succ)$. The associative product $\bar{\star}$ is then a cluster of four operations $x \bar{\star} y := x \nearrow  y+  x \searrow y+ x \swarrow y + x \nwarrow y $. 
In terms of categories, the following diagram commutes:
\begin{center}

$
\begin{array}{ccc}
\textsf{Quadri.} & \stackrel{F_v}{\longrightarrow} & \textsf{DiDend.} \\ 
 & & \\ 
 F_h \downarrow & \searrow f     & \downarrow F   \\
& & \\ 
\textsf{DiDend.} & \stackrel{F}{\longrightarrow}  & \textsf{As.}
\end{array} 
$
\end{center}
This diagram is commutative, i.e., $FF_v = f= FF_h$, where $F_v$ (resp. $F_h$) is the functor giving the vertical structure (resp. the horizontal structutre) of a quadri-algebra.
Section \ref{Octo} is devoted to
a generalisation of this structure obtained by considering three commuting Baxter operators.
\subsection{Construction of $\epsilon$-bialgebras from weighted directed graphs}
Let $X,V$ be two non-empty sets. Let $As(V)$ be the free associative algebra generated by $V$. By $As(V)\bra X \ket$, we mean the non-commutative algebra of polynomials with indeterminates belonging to $X$ with coefficients in $As(V)$, i.e., $ x \in As(V)\bra X \ket$ is a linear combinaison of the form $v_1v_2X_1v_3 \ldots v_p X_2 \ldots$, with $X_i \in X$ and $v_j \in As(V)$.
\begin{prop}
Let $X,V$ be two non-empty sets. Let $As(V)$ be  the free associative algebra generated by $V$. Suppose there exists a coassociative co-operation $\Delta: kX \xrightarrow{} kX^{\otimes 2}$. Then, the linear map
$\Delta_{\sharp}: As(V)\bra X \ket \xrightarrow{} As(V)\bra X \ket ^{\otimes 2}$ defined by $\Delta_{\sharp}(v)=0$ on $V$, by $\Delta_{\sharp}(X) := \Delta(X)$ on $X$ and extended by induction on $As(V)\bra X \ket$ by,
$$\Delta_{\sharp}(a_1a_2):= \Delta_{\sharp}(a_1)a_2 + a_1\Delta_{\sharp}(a_2), \ \forall a_1,a_2 \in As(V)\bra X \ket,$$
 is coassociative.
\end{prop}
\Proof
Straightforward.
\eproof
\Rk [\textbf{Interpretation}]
As a co-operation ``breaks down'' an element say $x$ into pieces $x_{(1)}$ and $x_{(2)}$ in several manners, we can say that elements of $V$ are indecomposable but not those of $X$. 

\noindent
Let $G$ be a weighted directed graph and $G_n$ be the set of path of length $n$, in particular $G_0$ is the vertex set  and $G_1$ the arrow set. A directed arrow $e_i \longrightarrow e_j$, with $e_i, \ e_j \in G_0$, will be denoted by $a_{(i,j)}$. A directed graph is said to be {\it{weighted}} if it is equipped with a weight map, i.e., with a map $w$ from $G_1$ to $k$.
The relations
$e_ie_j =1$ if $i=j$ and 0 otherwise, $a_{(i,j)} a_{(k,l)} = 0$ if $j \not= k$, $e_ia_{(k,l)} = a_{(k,l)}$ if $i=k$ and 0 otherwise and $a_{(k,l)}e_i = a_{(k,l)}$ if $i=l$ and 0 otherwise, for every $e_i, \ e_j \in G_0$ and $a_{(i,j)}, \ a_{(k,l)} \in G_1$, turn the space $kG := \oplus_{n=0} ^{ \infty} kG_n$ into a $k$-algebra called the path algebra. The product is then the concatenation of path whenever possible. We improve a construction given in \cite{Aguiar} in the case of directed graphs (trivially weighted). Define the following co-operation $\Delta$, for any $e_i \in G_0$ by $\Delta(e_i):=0$, for any weighted arrow $a_{(i,j)}$ by $\Delta a_{(i,j)}:= w(a_{(i,j)}) e_i \otimes e_j$ and for any weighted path $\alpha:= a_{(i_1,i_2)} a_{(i_2,i_3)} \ldots a_{(i_{n-1},i_n)}$ by, 
$$\Delta (\alpha):= \Delta (a_{(i_1,i_2)}) a_{(i_2,i_3)} \ldots a_{(i_{n-1},i_n)} + \ldots + a_{(i_1,i_2)} a_{(i_2,i_3)} \ldots a_{(i_{n-2},i_{n-1})}\Delta a_{(i_{n-1},i_n)}.$$

\begin{prop}[\cite{Lertribax}]
\label{dg}
Equipped with this co-operation, any weighted directed graph carries a $\epsilon$-bialgebra structure.
\end{prop}
\Rk
Elements which are unbreakable with regards to the co-operation are the vertices of the graph. 
\Rk \textbf{[Insertion of sequences and Virus attacks ]} \ \
From $\epsilon$-bialgebras, pre-Lie algebras can be constructed.
Let $(A, \ \mu, \ \Delta)$ be a $\epsilon$-bialgebra. Add a pre-Lie structure on $A$ via the following operation \cite{Aguiar}:
$$ \alpha \bowtie \alpha' := \alpha'_{(1)}\alpha \alpha'_{(2)}, \ \forall \alpha, \ \alpha' \in A, \ \textrm{with} \ \Delta(\alpha):= \alpha_{(1)} \otimes \alpha_{(2)}.  $$
Our structure on weighted directed graphs can be interpreted as follows.
When virus attack a cell, it breaks its D.N.A. chain $\alpha$ into two pieces $\alpha_{(1)}$ and $\alpha_{(2)}$ to insert its own D.N.A. chain. The result gives the D.N.A. chain $\alpha_{(1)}\alpha \alpha_{(2)}$. However, to do that, we can say that virus need an energy $\lambda_{12}$, or breaking the D.N.A. chain $\alpha$ of the host cell at this place can be done with a probability $\lambda_{12}$. 
\begin{center}
\includegraphics*[width=9cm]{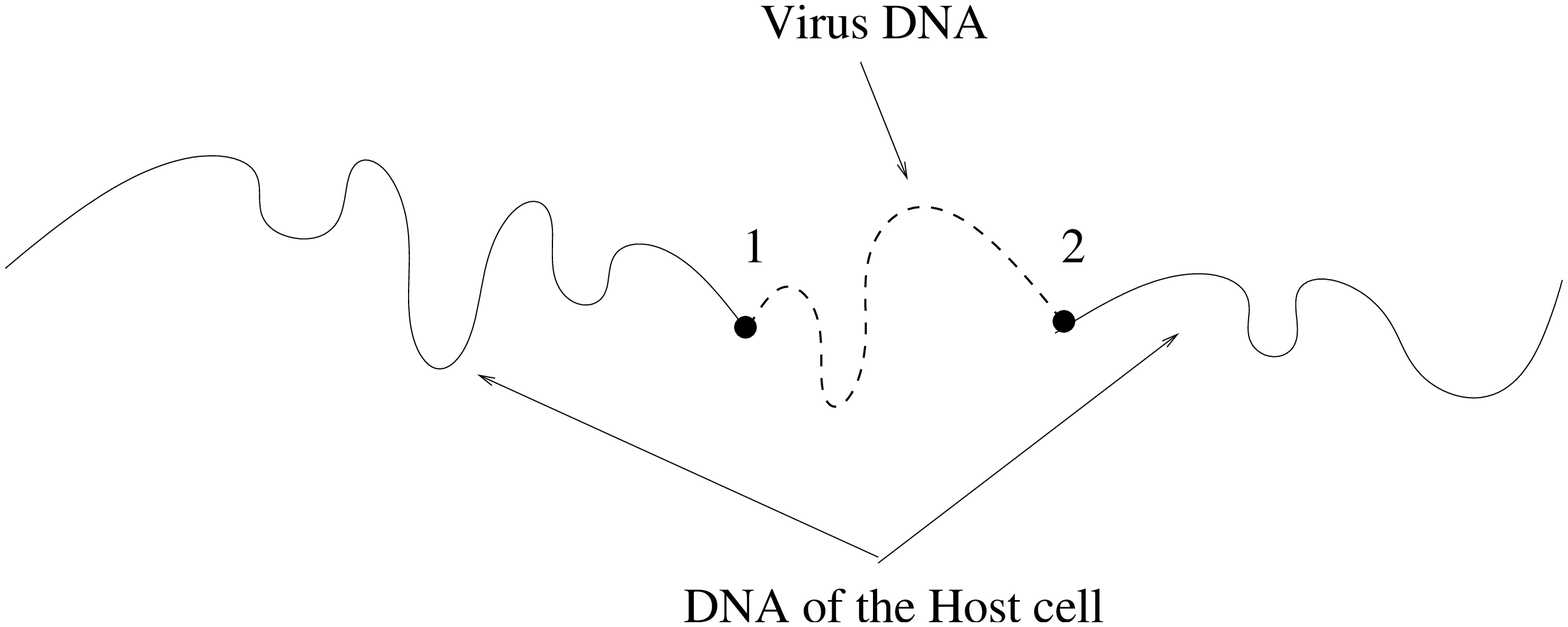}

\begin{scriptsize}
\textbf{
Insertion of a viral D.N.A. chain into the D.N.A. of a Host cell. When all possibilities are taken into account, this insertion mechanism has a pre-Lie structure coming from a natural $\epsilon$-bialgebra structure.}
\end{scriptsize}
\end{center}
With our modification of the coproduct introduced by Aguiar on graphs, such attacks and the probability or the energy to break the D.N.A. chain at this place rather than another place can be modelled. 
The coproduct,
viewed as a map which ``reads sequentially'' the D.N.A. chain,
models all the possibilities to open the D.N.A. chain
of the host cell with a given probability or a given energy. This entails that the insertion of a viral D.N.A. chain into the genetic code of the host cell has a natural pre-Lie structure.

This model can be improved by introducing more parameters to open a chain. This leads to the construction of $n$ coproducts $\Delta_k$, $k=1, \ldots n,$ such as the one defined in this example, but with the family of weights $(\lambda_{ij})_{(i \longrightarrow j \ \in G_1)}$ replaced by other families  $(\lambda^{(k)}_{ij})_{(i \longrightarrow j \ \in G_1)}$, where $k=1, \ldots, n$. Observe then that these coproducts obey
the following $n^2$ relations:
$$(\Delta_k \otimes id) \Delta_l := (id \otimes \Delta_l) \Delta_k,  $$
for all $k,l =1, \ldots n$. In \cite{Lertribax}, such structures have been used to produce formal deformation \`a la Gerstenhaber of quadri-algebras. See also Section~4. Such structures appear in \cite{Lertribax,codialg1} and in free probability \cite{Voi}. 
  
\subsection{Other constructions}
\label{2-7}
Consider the $k$-algebra $M_2(k)$ and the following co-operations.
\[\Delta_0 (\begin{pmatrix}
 a & b\\
c & d
\end{pmatrix}) = 
\begin{pmatrix}
 0 & a \\
 0 & c
\end{pmatrix} \otimes 
\begin{pmatrix}
 0 & 1 \\
 0 & 0
\end{pmatrix} -
\begin{pmatrix}
 0 & 1 \\
 0 & 0
\end{pmatrix} \otimes 
\begin{pmatrix}
 c & d \\
 0 & 0
\end{pmatrix}, \]
\[\Delta_1 (\begin{pmatrix}
 a & b\\
c & d
\end{pmatrix}) = 
\begin{pmatrix}
 b & 0 \\
 d & 0
\end{pmatrix} \otimes 
\begin{pmatrix}
 0 & 0 \\
 1 & 0
\end{pmatrix} -
\begin{pmatrix}
 0 & 0 \\
 1 & 0
\end{pmatrix} \otimes 
\begin{pmatrix}
 0 & 0 \\
 a & b
\end{pmatrix}, \]
Observe that $(M_2(k), \ \Delta_0)$, found in \cite{Aguiar}, and $(M_2(k), \ \Delta_1)$  are both $\epsilon$-bialgebras. The same can be obtained on $M_n(k)$ by defining $n$ coproducts $\Delta_{(a,b)}$ with $a+b =n+1$ such that,
$$ \Delta_{(a,b)} X := \sum_k X_{ka} E_{kb} \otimes E_{ab} - E_{ab} \otimes \sum_{k'} X_{bk'} E_{ak'},$$
where $X := X_{ij}E_{ij}$ with $X_{ij} \in k$ and $E_{ij} \in M_n(k)$ are the matrices with a 1 at the intersection of the row $i$ and the column $j$. For every $(a,b)$ with $a+b =n+1$, $(M_n(k), \ \Delta_{(a,b)})$ is a $\epsilon$-bialgebra.
\section{Algebraic structures related to two entangled coproducts.}
We continue our exploration of these particular bialgebras by considering 
this time $k$-vector spaces equipped with two co-operations. 
\subsection{$L$-Coalgebras $(A, \mu, \Delta_1, \Delta_2)$, with $(A, \mu,\Delta_2)$ a $\epsilon'[R]$-bialgebra.}
Let $(A, \ \mu)$ be an associative algebra equipped with two different co-operations $\Delta_1, \ \Delta_2: A \xrightarrow{} A^{\otimes 2} $. The $k$-algebra $\textsf{End}(A)$ has two convolution products denoted by $*_1$ and $*_2$. In the sequel, we will use four linear maps $\beta_1, \ \beta_2, \ \gamma_1, \ \gamma_2:\textsf{End}(A) \xrightarrow{} \textsf{End}(A)$ defined by $\beta_1(T):=id *_1 T$, $\beta_2(T):=id *_2 T$, $\gamma_1(T):= T *_1 id$, $\gamma_2(T):= T *_2 id,$ for all $T \in \textsf{End}(A)$.
\begin{prop}
\label{Prop 3.1}
Let $(A, \mu, \Delta_1, \Delta_2)$ be an associative algebra equipped with two co-operations $\Delta_1, \Delta_2$. Suppose $(A, \mu, \Delta_2)$ is a $\epsilon'[R]$-bialgebra. Then,
$$(Bax. 12) \ \ \ \ \beta_2(T)\beta_1(S)=
\beta_1(\beta_2(T)S), \ \forall T, S \in \textsf{End}(A).$$
In addition, the relations between operations labelled by 1 and 2 are stated in the following 2 by 2 matrix $(M_{4,\beta})_{ij}$. For all $T,R,S \in \textsf{End}(A)$, 
$$  R \overleftarrow{ \star_{\beta_1}} (S \prec_{\beta_2} T) = (R \overleftarrow{\star_{\beta_1}} S) \prec_{\beta_2} T,  \ \ \ 
R \overleftarrow{\star_{\beta_2}} (S \prec_{\beta_1} T) = (R \overleftarrow{\star_{\beta_2}} S ) \prec_{\beta_1} T,$$
$$ R \prec_{\beta_1} (S \overleftarrow{\star_{\beta_2}} T) = (R \prec_{\beta_2} S ) \prec_{\beta_1} T, \ \ \
R \overleftarrow{\star_{\beta_2}} (S \overleftarrow{\star_{\beta_1}} T) = (R \overleftarrow{\star_{\beta_2}} S) \overleftarrow{\star_{\beta_1}} T. $$
If $(A, \mu, \Delta_1)$ is also a $\epsilon'[R]$-bialgebra then
$(\textsf{End}(A), \ \overleftarrow{\star_{\beta_2}} \dashrightarrow \overleftarrow{\star_{\beta_1}})$ is an associative $L$-algebra.
\end{prop}
\Proof
Let $(A, \mu, \Delta_1, \Delta_2)$ be an algebra equipped with two co-operations $\Delta_1, \Delta_2$. Suppose $(A, \mu, \Delta_2)$ is a $\epsilon'[R]$-bialgebra. Fix $T,S \in \textsf{End}(A)$ and $x \in A$. Set $\Delta_1(x) := x_{(1)} \otimes x_{(2)}$ and $\Delta_2(x) := x^{(1)} \otimes x^{(2)}$.
On the one hand,
\begin{eqnarray*}
\beta_2(T)\beta_1(S)(x) &:=& \beta_2(T) (x_{(1)} S( x_{(2)} )), \\
&=& \mu (id \otimes T) \Delta_2 (x_{(1)} S( x_{(2)} )), \\
&=& \mu (id \otimes T) (x_{(1)} S( x_{(2)} )^{(1)} \otimes S( x_{(2)} )^{(2)} ),\\
&=& x_{(1)} S( x_{(2)} )^{(1)} T (S( x_{(2)} )^{(2)} ).
\end{eqnarray*}
On the other hand,
\begin{eqnarray*}
\beta_1(\beta_2(T)S)(x) &:=& \mu (id \otimes \beta_2(T)S) \Delta_1(x), \\
&=& x_{(1)} \beta_2(T)S(x_{(2)}), \\
&=& x_{(1)} \mu(id \otimes T)\Delta_2 (S(x_{(2)})), \\
&=& x_{(1)} \mu(id \otimes T) S(x_{(2)})^{(1)}  \otimes S(x_{(2)})^{(2)},\\
&=& x_{(1)}  S(x_{(2)})^{(1)} T( S(x_{(2)})^{(2)}).
\end{eqnarray*}
Therefore, $$(Bax. 12) \ \ \ \beta_2(T)\beta_1(S)=
\beta_1(\beta_2(T)S), \ \forall T, S \in \textsf{End}(A).$$
The other equalities are straightforward consequences of $(Bax. 12)$. If $(A, \mu, \Delta_1)$ is also a $\epsilon'[R]$-bialgebra then $(\textsf{End}(A), \ \overleftarrow{\star_{\beta_1}})$ is an associative algebra, see for instance Theorem \ref{ThAD}. As $\overleftarrow{\star_{\beta_2}}$ is also associative and is linked to $\overleftarrow{\star_{\beta_1}}$ via the equality $(M_{4, \beta})_{22}$, the $k$-vector space
$(\textsf{End}(A), \ \overleftarrow{\star_{\beta_2}} \dashrightarrow  \overleftarrow{\star_{\beta_1}})$ is an associative $L$-algebra.
\eproof
\begin{prop}
\label{Prop 3.2}
Let $(A, \mu, \Delta_1 \dashrightarrow \Delta_2)$ be a $L$-coalgebra, i.e., whose coproducts verify:
$$(\Delta_1 \otimes id)\Delta_2=(id \otimes \Delta_2)\Delta_1.$$  
Suppose $(A, \mu, \Delta_2)$ is a $\epsilon[R]$-bialgebra. Then,
$$(Bax. \ 3) \ \ \ \gamma_2(T)\gamma_1(S)=
\gamma_2(T\gamma_1(S)).$$
In addition, the operations labelled by 1 and 2 are stated in the following 2 by 2 matrix $(M_{4,\gamma})_{ij}$. For all $R, S, T \in \textsf{End}(A)$,
$$ R \succ_{\gamma_2} (S  \overrightarrow{\star_{\gamma_1}} T)= (R \succ_{\gamma_2} S) \overrightarrow{\star_{\gamma_1}} T, \ \ \ \
R \succ_{\gamma_1} (S \overrightarrow{\star_{\gamma_2}} T)= (R \succ_{\gamma_1} S) \overrightarrow{\star_{\gamma_2}} T,$$
$$R \overrightarrow{\star_{\gamma_2}} (S \overrightarrow{\star_{\gamma_1}} T) = (R \overrightarrow{\star_{\gamma_2}} S) \overrightarrow{\star_{\gamma_1}} T, \ \ \ \
R \succ_{\gamma_2}(S \succ_{\gamma_1} T) =(R \overrightarrow{\star_{\gamma_1}} S) \succ_{\gamma_2} T.$$
Moreover, for all $i,j= \{1,2 \}$ and $R, S, T \in \textsf{End}(A)$,
$$R \overleftarrow{ \star_{\beta_i} }
(S \overrightarrow{\star_{\gamma_j}} T)=
(R \overleftarrow{\star_{\beta_i}} S)
\overrightarrow{\star_{\gamma_j}} T$$  
$$R \succ_{\gamma_i} (S \prec_{\beta_j} T)=
(R \succ_{\gamma_i} S) \prec_{\beta_j} T.$$
If $(A, \mu, \Delta_1)$ is a $\epsilon[R]$-bialgebra then $(\textsf{End}(A), \ \overrightarrow{\star_{\gamma_2}} \dashrightarrow \overrightarrow{\star_{\gamma_1}})$ is an associative $L$-algebra.
\end{prop}
\Proof
Let $(A, \mu, \Delta_1 \dashrightarrow \Delta_2)$ be a $L$-coalgebra.
Suppose $(A, \mu, \Delta_2)$ is a $\epsilon[R]$-bialgebra. Fix $T,S \in \textsf{End}(A)$ and $x \in A$. Set $\Delta_1(x) := x_{(1)} \otimes x_{(2)}$ and $\Delta_2(x) := x^{(1)} \otimes x^{(2)}$. On the one hand,
\begin{eqnarray*}
\gamma_2(T)\gamma_1(S)(x) &=& \mu (T \otimes id)\Delta_2(S(x_{(1)})  x_{(2)}), \\
&=&\mu (T \otimes id) (S(x_{(1)}) (x_{(2)})^{(1)} \otimes (x_{(2)})^{(2)}), \\
&=& T(S(x_{(1)}) (x_{(2)})^{(1)}) (x_{(2)})^{(2)}.
\end{eqnarray*}
On the other hand,
\begin{eqnarray*}
\gamma_2(T\gamma_1(S))(x) &=& \mu(T\gamma_1(S) \otimes id)\Delta_2 (x), \\
&=& T\gamma_1(S)(x^{(1)})  x^{(2)}, \\
&=& T \mu(S \otimes id) (x^{(1)})_{(1)} \otimes (x^{(1)})_{(2)})  x^{(2)}, \\
&=& T(S  ((x^{(1)})_{(1)}) (x^{(1)})_{(2)} )       x^{(2)}, \\
\end{eqnarray*}
whence the equality $Bax. \ 3$ since $\Delta_1 \dashrightarrow \Delta_2$. The matrix $M_{4, \gamma}$ is obtained by using $Bax. \ 3$.  
If $(A, \mu, \Delta_1)$ is a $\epsilon[R]$-bialgebra then,
$(\textsf{End}(A), \ \overrightarrow{\star_{\gamma_1}})$ is an associative algebra. As $\overrightarrow{\star_{\gamma_1}}$ and $\overrightarrow{\star_{\gamma_2}}$ are related by the equation $(M_{4, \gamma})_{21}$, the $k$-vector space
$(\textsf{End}(A), \ \overrightarrow{\star_{\gamma_2}} \dashrightarrow \overrightarrow{\star_{\gamma_1}})$ is an associative $L$-algebra.
\eproof
\begin{prop}
\label{RRcat}
Let $(A, \mu, \Delta_1 \dashrightarrow \Delta_2)$ be a $L$-coalgebra. Suppose $(A, \mu, \Delta_1)$ is a $\epsilon'[R]$-bialgebra and $(A, \mu, \Delta_2)$ is a $\epsilon[R]$-bialgebra. Then, there exists a $( [M]_{ij})_{(i,j)}$-algebra structure on $\textsf{End}(A)$. 
\end{prop}
\Proof
As $(A, \mu, \Delta_1)$ is a $\epsilon'[R]$-bialgebra, the linear map $\beta_1$ is a right Baxter operator. As $(A, \mu, \Delta_2)$ is a $\epsilon[R]$-bialgebra the linear map $\gamma_2$ is a left Baxter operator. From $(id \otimes \Delta_2)\Delta_1=(\Delta_1 \otimes id)\Delta_2,$ we get the commutation relation
$\beta_1\gamma_2 = \gamma_2\beta_1.$
Apply Theorem \ref{RL Com} to conclude.
\eproof

\noindent
In terms of category, the following functor is obtained,
$$ [\epsilon'[R] \dashrightarrow \epsilon[R]]-\textsf{Bialg.} \longrightarrow (\textsf{[M]}_{(ij)}),$$
where $[\epsilon'[R] \dashrightarrow \epsilon[R]]-\textsf{Bialg.}$ denotes the category of those $L$-coalgebras verifying hypotheses of Proposition \ref{RRcat}.
\begin{exam}{}
We produce an example of $k$-vector space $(A, \mu, \ \Delta_1 \dashrightarrow \Delta_2)$ where both $(A, \mu, \ \Delta_1)$ and $(A, \mu, \ \Delta_2)$ are $\epsilon[R]$-bialgebras. 
Consider the Example \ref{$M_2(k)$}. Let $I_2$ be the identity matrix  and $X$ be any non-null $2 \times 2$-matrix. Then, the co-operation given by $\Delta_2 (M) := M \otimes X$ is coassociative and obey the identity:
$$ (\Delta_1 \otimes id)\Delta_2 := (id \otimes \Delta_2)\Delta_1. $$
The $k$-vector space $(M_2(k), \ \Delta_1 \dashrightarrow \Delta_2)$ is such that both $(M_2(k),  \ \Delta_1)$ and $(M_2(k), \ \Delta_2)$ are $\epsilon[R]$-bialgebras.
\end{exam}
\subsection{Dendriform dialgebras from $L$-coalgebras }
Let us construct dendriform algebras from non-coassociative co-operations. 
\begin{prop}
\label{epdend}
Let $(A, \ \mu, \ \Delta_1 \dashrightarrow \Delta_2)$ be a $L$-coalgebra, i.e., such that, $$ (\Delta_1 \otimes id)\Delta_2 = (id \otimes \Delta_2)\Delta_1.$$ Suppose $(A, \ \mu, \ \Delta_1)$ is a $\epsilon'[R]$-bialgebra and $(A, \ \mu, \ \Delta_2)$ is a $\epsilon$-bialgebra. Define the following binary operations $\prec, \ \succ : \textsf{End}(A)^{\otimes 2} \xrightarrow{} \textsf{End}(A)$ by,
$$T \prec S := \beta_1(T) \gamma_2(S) \ \ and \ \ T \succ S := \gamma_2\beta_1( T)S, \ \ \ \forall T,S \in \textsf{End}(A).$$
Then, $(\textsf{End}(A), \ \prec, \ \succ)$ is a dendriform dialgebra. 
\end{prop}
\Proof
As $(A, \ \mu, \ \Delta_1)$ is a $\epsilon'[R]$-bialgebra, the linear map $\beta_1$ is a right Baxter operator and $(\textsf{End}(A), \ \overleftarrow{*_{\beta_1}})$ is an associative algebra. 
As $(A, \ \mu, \ \Delta_2)$ is a $\epsilon$-bialgebra, $\gamma_2$
is a Baxter operator which commutes with $\beta_1$ since $(\Delta_1 \otimes id)\Delta_2 = (id \otimes \Delta_2)\Delta_1.$
This entails that $\gamma_2$ is a Baxter operator on the associative algebra $(\textsf{End}(A), \ \overleftarrow{*_{\beta_1}})$. Therefore, $(\textsf{End}(A), \ \prec, \ \succ)$ is a dendriform dialgebra. 
\eproof

\noindent
In terms of category, the following functor is obtained,
$$ [\epsilon'[R] \dashrightarrow \epsilon]-\textsf{Bialg.} \longrightarrow \textsf{DiDend.},$$
where $[\epsilon'[R] \dashrightarrow \epsilon]-\textsf{Bialg.}$ denotes the category of those $L$-coalgebras verifying hypotheses of Proposition \ref{epdend}.
\begin{prop}
\label{Dendeps}
Let $(A, \ \mu, \ \Delta_1 \dashrightarrow \Delta_2)$ be a $L$-coalgebra i.e., such that, $$ (\Delta_1 \otimes id)\Delta_2 = (id \otimes \Delta_2)\Delta_1.$$ Suppose $(A, \ \mu, \ \Delta_1)$ is a $\epsilon$-bialgebra and $(A, \ \mu, \ \Delta_2)$ is a $\epsilon[R]$-bialgebra. Define the following binary operations $\succ, \ \prec: \textsf{End}(A)^{\otimes 2} \xrightarrow{} \textsf{End}(A)$ by,
$$T \succ S := \beta_1(T) \gamma_2(S) \ \ and \ \ T \prec S := T \beta_1(\gamma_2(S)), \ \ \ \forall T,S \in \textsf{End}(A).$$
Then, $(\textsf{End}(A), \ \prec, \ \succ)$ is a dendriform dialgebra. 
\end{prop}
\Proof
As $(A, \ \mu, \ \Delta_1)$ is a $\epsilon$-bialgebra, the linear map $\beta_1$ is a Baxter operator on $\textsf{End}(A)$. 
As $(A, \ \mu, \ \Delta_2)$ is a $\epsilon[R]$-bialgebra, the linear map $\gamma_2$
is a left Baxter operator which commutes with $\beta_1$ since $(\Delta_1 \otimes id)\Delta_2 = (id \otimes \Delta_2)\Delta_1.$
This entails that $\beta_1$ is a Baxter operator on the associative algebra $(\textsf{End}(A), \ \overrightarrow{*_{\gamma_2}})$. Therefore, $(\textsf{End}(A), \ \prec, \ \succ)$ is a dendriform dialgebra. 
\eproof

\noindent
In terms of category, the following functor is obtained,
$$ [\epsilon \dashrightarrow \epsilon[R]]-\textsf{Bialg.} \longrightarrow \textsf{DiDend.},$$
where $[\epsilon \dashrightarrow \epsilon[R]]-\textsf{Bialg.}$ denotes the category of those $L$-coalgebras verifying hypotheses of Proposition \ref{Dendeps}.
\begin{exam}{[Dendriform dialgebras from weighted directed graphs and research algorithms]}
Keep notation of Proposition \ref{dg}.
Recall the co-operation $\Delta$, for any $e_i \in G_0$ is defined by $\Delta(e_i):=0$, for any weighted arrow $a_{(i,j)}$ by $\Delta a_{(i,j)}:= w(a_{(i,j)}) e_i \otimes e_j$ and for any weighted path $\alpha:= a_{(i_1,i_2)} a_{(i_2,i_3)} \ldots a_{(i_{n-1},i_n)}$ by, 
$$\Delta (\alpha):= \Delta (a_{(i_1,i_2)}) a_{(i_2,i_3)} \ldots a_{(i_{n-1},i_n)} + \ldots + a_{(i_1,i_2)} a_{(i_2,i_3)} \ldots a_{(i_{n-2},i_{n-1})}\Delta a_{(i_{n-1},i_n)}.$$
Equipped with this co-operation, any weighted directed graph carries a $\epsilon$-bialgebra structure. 
For any path $\alpha$ define the coproduct $\Delta_M$ by $\Delta_M(\alpha):= \alpha \otimes \sum_{ij} a_{ij}$ for any path $\alpha$. Equipped with this co-operation, any weighted directed graph carries a $\epsilon[R]$-bialgebra structure. Observe
that,
$$(\Delta \otimes id)\Delta_M = (id \otimes \Delta_M) \Delta. $$

The map $(\mu \otimes id)(id \otimes \Delta)\Delta_M$ can be interpreted as a research map throughout the graph $G$ since for any path $\alpha$, $(\mu \otimes id)(id \otimes \Delta)\Delta_M(\alpha):= \sum_j \lambda_{t(\alpha)j} \alpha \otimes e_{j}$, with $\lambda_{t(\alpha)j}$ a non-null scalar if there is an arrow $t(\alpha) \xrightarrow{} j $ with weight $\lambda_{t(\alpha)j}$ belonging to $G_1$.
Consider the following weighted directed graph $G$, with weights $\lambda_{ij}$ different from zero if $i \longrightarrow j \in G_1$ and zero otherwise. 
\begin{center}
\includegraphics*[width=6cm]{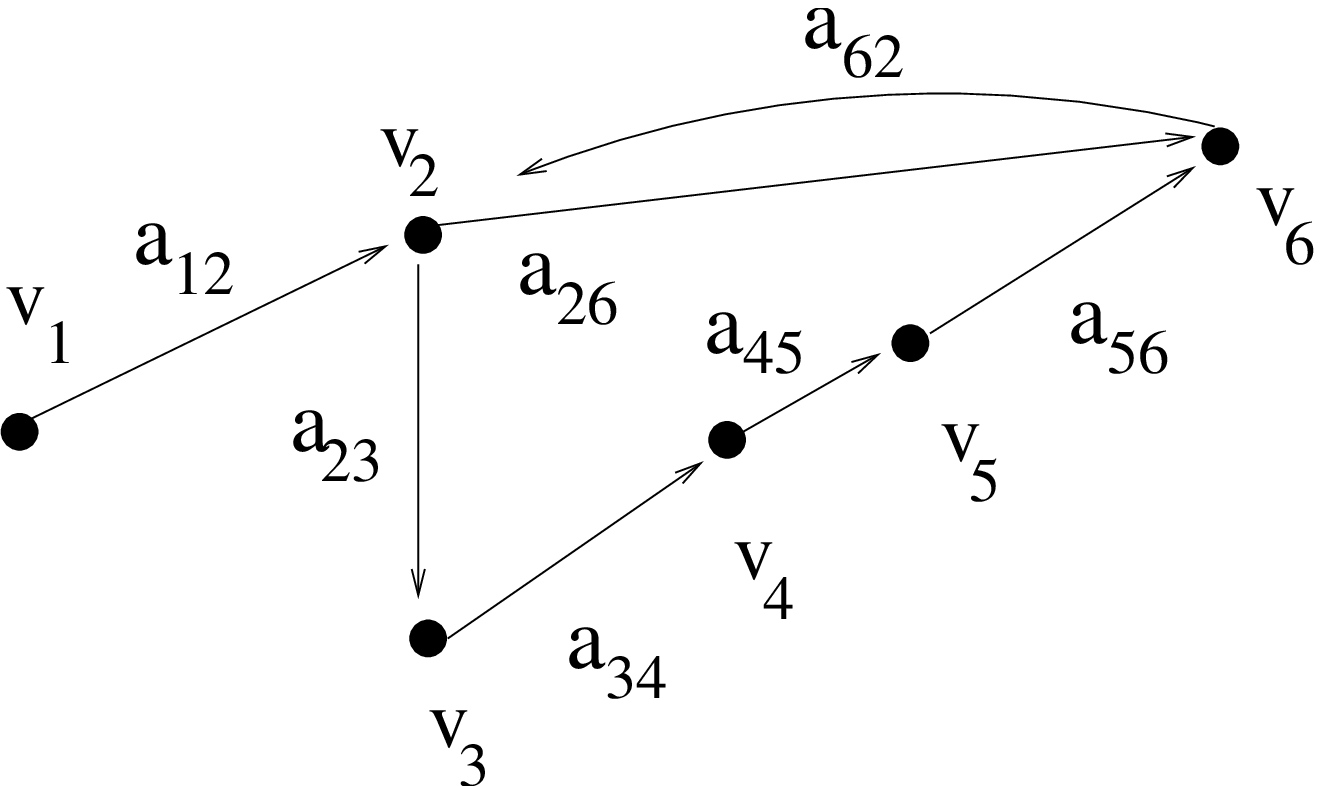}

\begin{scriptsize}
\textbf{The weighted directed graph $G$ (no weight represented).}
\end{scriptsize}
\end{center}
The operator $(\mu \otimes id)(id \otimes \Delta)\Delta_M$ applied to 
$\alpha := a_{12}a_{23}a_{34}$ gives,
$$ a_{12}a_{23}a_{34} \xrightarrow{\Delta_M} \sum_{i,j=1}^6 \ a_{12}a_{23}a_{34} \otimes a_{ij} \xrightarrow{id \otimes \Delta} \sum_{i,j=1}^6 \ \lambda_{ij} \ a_{12}a_{23}a_{34} \otimes e_i \otimes e_j
\xrightarrow{\mu \otimes id} \lambda_{45} \ a_{12}a_{23}a_{34} \otimes e_5.$$

\end{exam}
\subsection{Entanglement of left and right bialgebras}
\begin{prop}
\label{LRbb}
Let $(A, \ \mu, \Delta_1 \dashrightarrow \Delta_2)$ be a $L$-coalgebra such that $(A, \ \mu, \Delta_1)$ is a $\epsilon'[R]$-bialgebra and $(A, \ \mu, \Delta_2)$ is a $\epsilon'[L]$-bialgebra. Then, there exists a $([M_1]_{ij})_{(i,j)}$-algebra structure on $\textsf{End}(A)$.
\end{prop}
\Proof
As $(A, \ \mu, \Delta_2)$ is a $\epsilon'[L]$-bialgebra, $\gamma_2$ is a right Baxter operator. So is
 $\beta_1$ since $(A, \ \mu, \Delta_1)$ is a $\epsilon'[R]$-bialgebra. The two right Baxter operators commute, i.e., $\gamma_2\beta_1 = \beta_1\gamma_2$, since the notation $\Delta_1 \dashrightarrow \Delta_2$ means $(\Delta_1 \otimes id) \Delta_2 = (id \otimes \Delta_2) \Delta_1$. Consequently, $\gamma_2$ is a right Baxter operator on the $L$-anti-dipterous
$(\textsf{End}(A), \ \overleftarrow{*}_{\beta_1}, \prec_ {\beta_1})$. Corollary \ref{Transp RR} claims that $\textsf{End}(A)$ has a structure of $([M_1]_{ij})_{(i,j)}$-algebra.
\eproof

\noindent
In terms of category, the following functor is obtained,
$$ [\epsilon'[R] \dashrightarrow \epsilon'[L]]-\textsf{Bialg.} \longrightarrow [\textsf{M}_{1}]_{(ij)},$$
where $[\epsilon'[R] \dashrightarrow \epsilon'[L]]-\textsf{Bialg.}$ denotes the category of those $L$-coalgebras verifying hypotheses of Proposition \ref{LRbb}.

\begin{prop}
\label{LRb}
Let $(A, \ \mu, \Delta_1 \dashrightarrow \Delta_2)$ be a $L$-coalgebra such that $(A, \ \mu, \Delta_1)$ is a $\epsilon[L]$-bialgebra and $(A, \ \mu, \Delta_2)$ a $\epsilon[R]$-bialgebra. Then, there exists a $([M_2]_{ij})_{(i,j)}$-algebra structure on $\textsf{End}(A)$.
\end{prop}
\Proof
As $(A, \ \mu, \Delta_2)$ is a $\epsilon[R]$-bialgebra, $\gamma_2$ is a left Baxter operator. So is
 $\beta_1$ since $(A, \ \mu, \Delta_1)$ is a $\epsilon[L]$-bialgebra. The two left Baxter operators commute, i.e., $\gamma_2\beta_1 = \beta_1\gamma_2$, since the notation $\Delta_1 \dashrightarrow \Delta_2$ means $(\Delta_1 \otimes id) \Delta_2 = (id \otimes \Delta_2) \Delta_1$. Consequently, $\gamma_2$ is a left Baxter operator on the $L$-dipterous
$(\textsf{End}(A), \ \overrightarrow{*}_{\beta_1}, \succ_ {\beta_1})$. Proposition \ref{Prop Dip} claims that $\textsf{End}(A)$ has a structure of $([M_2]_{ij})_{(i,j=1,\ldots 3)}$-algebra.
\eproof

\noindent
In terms of category, the following functor is obtained,
$$ [\epsilon[L] \dashrightarrow \epsilon[R]]-\textsf{Bialg.} \longrightarrow [\textsf{M}_2]_{(ij)},$$
where $[\epsilon[R] \dashrightarrow \epsilon[L]]-\textsf{Bialg.}$ denotes the category of those $L$-coalgebras verifying hypotheses of Proposition \ref{LRb}.
\section{Splitting associative clusters into associative operations}
Let $(X,  \bullet)$ be a $k$-algebra and $(\bullet_i)_{1 \leq i \leq N}: X^{\otimes 2} \xrightarrow{} X$ be a family of binary operations on $X$. Recall that
$\bullet \longrightarrow \ \sum_i \bullet_i$ means $x \bullet y =  \sum_i x \bullet_i y$, for all $x,y \in X$, i.e.,
$\bullet$ {\it{splits}} into the $N$ operations $\bullet_1, \ldots, \bullet_N$ or is a {\it{cluster of $N$ (binary) operations}}. Such a cluster will be called associative if $\bullet$ is associative. For instance dendriform dialgebras have their associative operations which split into two operations $\star \longrightarrow  \ \prec + \succ$. Here, we will be interested in associative clusters which split into $n$ associative products.
\subsection{$n$-Hypercubic $\epsilon[R]$-bialgebras}
\label{hhhyyp}
\begin{defi}{[Hypercubic algebras]}
Fix a non-null integer $n$.
A {\it{$n$-hypercubic algebra}} $(A, \ \star_1, \ldots, \ \star_n)$ is a $k$-vector space equipped with $n$ binary operations
$\star_1, \ldots, \ \star_n: A^{\otimes 2} \xrightarrow{} A$, such that for all $x,y,z \in A$ and for all $i,j \in \{1, \ldots, n \}$:
$$ (Hyp.) \ \ \ (x \star_i y)\star_j z = x \star_i (y \star_j z). $$
Similarly,
a {\it{$n$-hypercubic coalgebra}} $(A, \ \Delta_1, \ldots, \ \Delta_n)$ is a $k$-vector space equipped with $n$ co-operations
$\Delta_1, \ldots, \ \Delta_n: A \xrightarrow{} A^{\otimes 2}$, such that for all $i,j \in \{1, \ldots, n \}$:
$$(CoHyp.) \ \ \ ( \Delta_i \otimes id)\Delta_j  = (id \otimes \Delta_j)\Delta_i. $$
A {\it{$n$-hypercubic $\epsilon[R]$-bialgebra}} $(A, \ \mu, \ \Delta_1, \ldots, \ \Delta_n)$ is a $n$-hypercubic coalgebra
$(A, \ \Delta_1, \ldots, \ \Delta_n)$ together with an associative algebra $(A, \ \mu)$ such that for all $i=1, \ldots,n,$ $(A, \ \mu, \ \Delta_i)$ are $\epsilon[R]$-bialgebras. Similarly for the left version.
A {\it{$n$-hypercubic $\epsilon$-bialgebra}} is a $n$-hypercubic coalgebra
$(A, \ \Delta_1, \ldots, \ \Delta_n)$ together with an associative algebra $(A, \ \mu)$ such that for all $i=1, \ldots,n,$ $(A, \ \mu, \Delta_i)$ are $\epsilon$-bialgebras. Such bialgebras have been used in \cite{Lertribax} to produce formal deformations of associative algebras such as dendriform di- and tri-algebras, quadri-algebras and ennea-algebras.
\end{defi}
\Rk
Observe that axioms $(Hyp.)$ entails the associativity of all the involved operations. Moreover these axioms are invariant by action of $M_n(k)$. More precisely, if $M \in M_n(k)$ and $X :=(\star_1, \ldots, \ \star_n)$ then $^t(M(^tX)) :=(\star'_1, \ldots, \ \star'_n)$ will still yield
associative products $\star'_1, \ldots, \ \star'_n$ which will satisfy the axioms stated above. In particular, $\sum_{i = 1}^n \star_i$ is still associative. 
If $(A, \ \star_1, \ldots, \ \star_n)$ is a $n$-hypercubic 
algebra, then $\star_j$ is a Hochschild 2-cocycle on the associative algebra $(A, \ \star_i )$, for all $i,j=1, \ldots n$. In addition, there also exists a notion of opposite structure. Indeed, define new operations $x \star^{op} _i y := y \star_i x$, for all $x,y \in A$ and $i=1, \ldots, n$ and observe that
$(A, \ \star^{op} _1, \ldots, \ \star^{op} _n)$ is still an hypercubic algebra.
\Rk
Since the operad associated with hypercubic algebras is binary, quadratic and non-symmetric,
the free $n$-hypercubic algebra on a $k$-vector space $V$ is given by:
$$ H^n(V):= \bigoplus_{p \geq 1} \ (H^n)_p \otimes V^{\otimes p}.$$
On one generator, one see that $H^n(k):= \bigoplus_{p \geq 1} \ (H^n)_p$. The dimension of $(H^n)_p$ is $n^p$. 
These operads are conjectured to be Koszul self-dual. This has already been proved in the case $p=2$ \cite{Richter} and $p=3$ \cite{LodayRonco}.

\noindent
Let us construct hypercubic algebras from non-coassociative coproducts.
\begin{theo}
Let $(A, \ \mu, \ \Delta_1, \dots, \Delta_n)$ be a $k$-vector space such that for all $i=1, \ldots, n$, $(A, \ \mu, \ \Delta_i)$ are $\epsilon'[R]$-bialgebras. Then, $(\textsf{End}(A), \ \overleftarrow{*}_{\beta_1}, \ldots, \ \overleftarrow{*}_{\beta_n})$ is a $n$-hypercubic algebra.
\end{theo} 
\Proof
As for all $i=1, \ldots, n$, $(A, \ \mu, \ \Delta_i)$ are $\epsilon'[R]$-bialgebras, $(\textsf{End}(A), \ \overleftarrow{*}_{\beta_i})$ are associative algebras. Apply now the last claim of Proposition \ref{Prop 3.1} to any pair $((A, \ \mu, \ \Delta_i)$, $(A, \ \mu, \ \Delta_j))$ to prove that $(\textsf{End}(A), \ \overleftarrow{*}_{\beta_1}, \ldots, \ \overleftarrow{*}_{\beta_n})$ is a $n$-hypercubic algebra.
\eproof
\begin{prop}
Let $(A, \ \mu, \ \Delta_1, \dots, \Delta_n)$ be  a $n$-hypercubic $\epsilon[R]$-bialgebra. Then, the $k$-vector space $(\textsf{End}(A), \ \overrightarrow{*}_{\gamma_1}, \ldots, \ \overrightarrow{*}_{\gamma_n})$ is a $n$-hypercubic algebra.
\end{prop} 
\Proof
Apply now the last claim of Proposition \ref{Prop 3.2} to any pair $((A, \ \mu, \ \Delta_i)$, $(A, \ \mu, \ \Delta_j))$.
\eproof
\begin{prop}
Let $(A, \ \mu, \ \Delta_1, \Delta_2)$ be a $2$-hypercubic $\epsilon$-bialgebra.
Then, the associative products $\overleftarrow{*}_{\beta_1}$ and $\overleftarrow{*}_{\beta_2}$ defined on $\textsf{End}(A)$ verify the relation:
$$(Aas.): \ \ \ \  (R *_{\beta_1} S)*_{\beta_2} T+
(R *_{\beta_2} S)*_{\beta_1} T =
R *_{\beta_1}(S *_{\beta_2} T) +
R *_{\beta_2}(S *_{\beta_1} T).$$
Moreover the operations $\star_{[\lambda_1, \ \lambda_2]} := \lambda_1(*_{\beta_1}) + \lambda_2(*_{\beta_2}): \textsf{End}(A) \xrightarrow{} \textsf{End}(A)$, where $\lambda_1, \ \lambda_2 \in k$ are still associative. 
\end{prop}
\Proof
Let $(A, \ \mu, \ \Delta_1, \Delta_2)$ be a $2$-hypercubic $\epsilon$-bialgebra. Fix $R,S,T  \in \textsf{End}(A)$. 
Observe that, 
$$\beta_1(T)\beta_2(S):= \beta_1(T\beta_2(S)) + \beta_2(\beta_1(T)S),$$
if $(\Delta_1 \otimes id)\Delta_2 = (id \otimes \Delta_2)\Delta_1$.
This gives $(Aas'.)$:
$$ (R *_{\beta_1} S) *_{\beta_2} T - R *_{\beta_1} (S *_{\beta_2 T}) =
\beta_2(R \beta_1(S))T + R \beta_2(\beta_1(S)T) - R \beta_1(\beta_2(S)T) - \beta_1(R \beta_2(S))T.$$
The equation $(Aas.)$ holds by noticing that $(Aas'.)$ is anti-symmetric with regards to the permutation $(2,1)$. The last assertion is straightforward.
\eproof
\Rk
The condition $(Aas.)$ claims that $*_{\beta_1}$ is a Hochschild 2-cocycle on the associative algebra $(\textsf{End}(A), \ *_{ \beta_2})$ or that
$*_{\beta_2}$ is a Hochschild 2-cocycle on the associative algebra $(\textsf{End}(A), \ *_{ \beta_1})$.
\Rk
Observe here that associative products $\star_{[\lambda_1, \ \lambda_2]} \rightarrow \lambda_1(\overleftarrow{*}_{\beta_1}) + \lambda_2(\overleftarrow{*}_{\beta_2})$ split into two associative ones. The following corollary improves this result.
\begin{coro}
Let $(A, \ \mu, \ \Delta_1, \ldots,  \ \Delta_n)$ be a $n$-hypercubic $\epsilon$-bialgebra.
Then, the associative products $\star := \sum_i  \lambda_i \overleftarrow{*}_{\beta_i}, \ \hat{\star} := \sum_i  \lambda'_i \overrightarrow{*}_{\gamma_i}: \textsf{End}(A) \xrightarrow{} \textsf{End}(A)$ are associative.
\end{coro}
\Proof
Straightforward. For instance, 
fix $R,S,T  \in \textsf{End}(A)$. 
Observe that, 
$$\gamma_i(T)\gamma_j(S):= \gamma_i(T\gamma_j(S)) + \gamma_j(\gamma_i(T)S),$$
if $(\Delta_j \otimes id)\Delta_i = (id \otimes \Delta_i)\Delta_j$. With
$\beta$ replaced by $\gamma$, the relation $(Aas.)$ still holds. Therefore $\hat{\star} := \sum_i  \lambda'_i \overrightarrow{*}_{\gamma_i}$ is associative.
\eproof
\Rk
Solving equation $(Aas.)$ entails the study of at least two types of algebras, the $n$-hypercubic algebras we have just seen and the $n$-circular algebras. Let us see the interest of such algebras and their co-versions.
\subsection{Circular algebras}
\begin{defi}{[circular algebras]}
A {\it{circular algebra}} $(A, \ \star_1, \ldots, \ \star_n)$ is a $k$-vector space equipped with $n$ binary operations
$\star_1, \ldots, \ \star_n: A^{\otimes 2} \xrightarrow{} A$, such that for all $x,y,z \in A$ and for all $i,j \in \{1, \ldots, n \}$:
$$ (x \star_i y)\star_j z = x \star_j (y \star_i z). $$
\end{defi}
\Rk
Observe that such a definition entails the associativity of all the involved operations and that axioms of a circular algebra are invariant by action of $M_n(k)$. More precisely, if $M \in M_n(k)$ and $X :=(\star_1, \ldots, \ \star_n)$ then $^t(M(^tX)) :=(\star'_1, \ldots, \ \star'_n)$ will still yield
associative products $\star'_1, \ldots, \ \star'_n$ which will satisfy the axioms stated above. In particular, $\sum_{i = 1}^n \star_i$ is still associative.
\Rk
If $(A, \ \star_1, \ldots, \ \star_n)$ is a $n$-circular
algebra, then $\star_j$ is a Hochschild 2-cocycle on the associative algebra $(A, \ \star_i )$.
\Rk
If there exists a common unit to all these products $\star_1, \ldots, \ \star_n$, then they are equal pairwise. 
\Rk
Since the non-symmetric operad associated with circular algebras is binary and quadratic,
the free circular algebra on a $k$-vector space $V$ is given by:
$$ C^n(V):= \bigoplus_{p \geq 1} \ (C^n)_p \otimes V^{\otimes p}.$$
On one generator, one see that $C^n(k):= \bigoplus_{p \geq 1} \ (C^n)_p$. The dimension of $C_p$ is $n^p$.
\subsubsection{Constructions of $n$-circular algebras via right Baxter operators}
\begin{lemm}
\label{lemm circ}
Let $(A, \mu)$ be an associative algebra.
Suppose there exist $n$ linear maps $\phi_1, \dots, \phi_n: A \xrightarrow{} A$ such that for all $x,y \in A$ and for $i,j \in \{1, \ldots n, \}$, 
$$ \phi_i(x) \phi_j(y) := \phi_i(\phi_j(x)y),$$ 
hold.
Define the binary operations $\star_1, \dots, \star_n$ by $x \star_1 y := \phi_1(x)y, \ldots, \ x \star_n y := \phi_n(x)y$.
Then, the $k$-vector space $(A, \ \star_1, \dots, \ \star_n)$ is a $n$-circular algebra. Similarly, Suppose these $n$-linear maps verify,
$$ \phi_i(x) \phi_j(y) := \phi_j(x\phi_i(y)),$$ 
for $i,j \in \{1, \ldots n, \}$. Then, the binary operations $\star'_1, \dots, \star'_n$ define by $x \star'_1 y := x\phi_1(y), \ldots, \ x \star'_n y := x\phi_n(y)$ turn
the $k$-vector space $(A, \ \star'_1, \dots, \ \star'_n)$ into a $n$-circular algebra.
\end{lemm}
\Proof
Straightforward.
\eproof
\begin{prop}
Let $(A, \ \mu, \ \Delta_1, \dots, \Delta_n)$ be a $n$-circular $\epsilon[R]$-bialgebra. Then, there exists a $n$-circular algebraic structure on $\textsf{End}(A)$.
\end{prop} 
\Proof
Recall that $\gamma_i(T):= T *_i id$, for all $T \in \ \textsf{End}(A)$ and $i=1, \ldots, n$. Fix $T, S \in \textsf{End}(A)$ and $a \in A$. Set $\Delta_1(a):= a_{(1)} \otimes a_{(2)}$ and $\Delta_1(a):= a^{(1)} \otimes a^{(2)}$. On the one hand, 
$
\gamma_2(T)\gamma_1(S)(a):= T(S(a_{(1)})(a_{(2)})^{(1)})(a_{(2)})^{(2)}$. On the other hand,
$\gamma_1(T\gamma_2(S))(a):= T \mu(S \otimes id)\Delta_2(a_{(1)})a_{(2)} = T(S((a_{(1)})^{(1)})(a_{(1)})^{(2)})a_{(2)}.$ The two relations are equal if $\Delta_2(a_{(1)}) \otimes a_{(2)} = a_{(1)} \otimes \Delta_2(a_{(2)})$. We obtain a family of operators $\gamma_1, \ldots, \gamma_n,$ verifying for all $i,j =1, \ldots,n,$ and for all $T,S \in \textsf{End}(A)$,
$$ \gamma_i(T)\gamma_j(S) = \gamma_j(T\gamma_i(S)).$$
Apply now Lemma \ref{lemm circ} to conclude.
\eproof
\subsubsection{Constructions of circular algebras}
We give two examples of $2$-circular algebras constructed from convolution products on the $k$-algebra $M_2(k)$. Fix $\lambda, \nu \in k$ and consider the following co-operations.
\[\Delta_1 (\begin{pmatrix}
 a & b\\
c & d
\end{pmatrix}) = 
\begin{pmatrix}
 0 & a \\
 0 & c
\end{pmatrix} \otimes 
\begin{pmatrix}
 0 & 1 \\
 0 & 0
\end{pmatrix} -
\begin{pmatrix}
 0 & 1 \\
 0 & 0
\end{pmatrix} \otimes 
\begin{pmatrix}
 c & d \\
 0 & 0
\end{pmatrix}, \]
\[
\Delta_2 (\begin{pmatrix}
 a & b\\
c & d
\end{pmatrix}) = 
\begin{pmatrix}
 0 & \lambda a \\
 0 & \nu c
\end{pmatrix} \otimes 
\begin{pmatrix}
 0 & 1 \\
 0 & 0
\end{pmatrix} -
\begin{pmatrix}
 0 & 1 \\
 0 & 0
\end{pmatrix} \otimes 
\begin{pmatrix}
 \nu c & \lambda d \\
 0 & 0
\end{pmatrix}.
\]
Then, $(\Delta_i \otimes id) \Delta_j = (id \otimes \Delta_i) \Delta_j$ holds for all $i,j:= 1,2$. Similarly,
\[\Delta'_1 (\begin{pmatrix}
 a & b\\
c & d
\end{pmatrix}) = 
\begin{pmatrix}
 b & 0 \\
 d & 0
\end{pmatrix} \otimes 
\begin{pmatrix}
 0 & 0 \\
 1 & 0
\end{pmatrix} -
\begin{pmatrix}
 0 & 0 \\
 1 & 0
\end{pmatrix} \otimes 
\begin{pmatrix}
 0 & 0 \\
 a & b
\end{pmatrix}, \]
\[
\Delta'_2 (\begin{pmatrix}
 a & b\\
c & d
\end{pmatrix}) = 
\begin{pmatrix}
 \lambda b &  0 \\
 \nu d & 0
\end{pmatrix} \otimes 
\begin{pmatrix}
 0 & 0 \\
 1 & 0
\end{pmatrix} -
\begin{pmatrix}
 0 & 0 \\
 1 & 0
\end{pmatrix} \otimes 
\begin{pmatrix}
  0 & 0 \\
 \nu a & \lambda b
\end{pmatrix}.
\]
The equalities $(\Delta'_i \otimes id) \Delta'_j = (id \otimes \Delta'_i) \Delta'_j$ hold for all $i,j := 1,2$.
The $k$-algebra $M_2(k)$ equipped with the two convolution products $*_1$ and $*_2$ or $*'_1$ and $*'_2$ is a 2-circular algebra. We do not know how to construct 2-circular $\epsilon$-bialgebras or their left and right versions. However, observe that $(M_2(k), \ \Delta_1)$ is a $\epsilon$-bialgebra \cite{Aguiar}. See also Subsection \ref{2-7}.
\section{Octo-algebras}
\label{Octo}
Up to now, we have handled two commuting Baxter operators. What kind of structures emerge if we work with three pairwise commuting Baxter operators? 
\subsection{Three different ways to obtain quadri-algebras from an associative cluster}
We weak hypotheses of Section~4 and introduce the $L$-circularity.
Let $V$ be a $k$-vector space equipped with two co-operations $\Delta_1$ and $\Delta_2$. The co-operation $\Delta_1$ is said to be {\it{$L$-circular}} with regards to $\Delta_2$ if:
\begin{enumerate}
\item {$\Delta_1 \dashrightarrow \Delta_2,$}
\item{$ (\Delta_1 \otimes id)\Delta_2 := (id \otimes \Delta_1)\Delta_2$.}
\end{enumerate} 
This condition will be denoted by the symbol $\Delta_1 \curvearrowright \Delta_2$. One of the interest of $L$-circularity is to produce three Baxter operators which commute pairwise in the case of a $L$-circular $\epsilon$-bialgebra, which is simply a $k$-vector space $(A, \ \mu, \Delta_1 \curvearrowright \Delta_2)$ where both $(A, \ \mu, \Delta_1)$ and $(A, \ \mu, \ \Delta_2)$ are $\epsilon$-bialgebras.
This gives birth to octo-algebras, which are simply associative algebras whose associative products are clusters of 8 operations related by 27 conditions. From such an octo-algebra, 3 quadri-algebras or 6 dendriform dialgebras sharing the same associative operation are constructed. Let us see how it works.
\begin{defi}{[Octo-algebra]}
\label{defocto}
An {\it{octo-algebra}} $O$ is a $k$-vector space equipped with 8 binary operations: $$\swarrow_1, \ \swarrow_2, \ \nwarrow_1, \ \nwarrow_2, \ \searrow_1, \ \searrow_2, \ \nearrow_1, \ \nearrow_2: O^{\otimes 2} \xrightarrow{} O.$$
To ease notation, the following sums are introduced.
\begin{eqnarray*}
x \prec_i y &:=& x \swarrow_i y + x \nwarrow_i y; \ \ \ \ \ x \bigwedge y:=x  \wedge_1 y+ x \wedge_2 y; \\
x \succ_i y &:=& x \searrow_i y + x \nearrow_i y;  \ \ \ \ \  x \bigvee y:= x \vee_1 y+ x \vee_2 y, \\                    
x \wedge_i y &:=& x \nwarrow_i y + x \nearrow_i y; \ \ \ \ \ x \ll y:= x \prec_1 y + x \prec_2 y, \\
x \vee_i y &:=& x \swarrow_i y+ x \searrow_i y; \ \ \ \ \ x \gg y:= x \succ_1 y + x \succ_2 y,  \\
x \circ_{12} y &:=& x \circ_{1} y + x \circ_{2} y, \ \ \ \ \ \ \ \ \textrm{with} \ \circ \in \{ \  \nearrow, \ \searrow, \ \swarrow, \ \nwarrow \},\\
x \ \Sigma_1 \ y &:=& x \vee_1 y + x \wedge_1 y = x \prec_1 y+ x \succ_1 y; \ \ \ \ \
x \ \Sigma_2 \ y := x \vee_2 y + x \wedge_2 y= x \prec_2 y+ x \succ_2 y; 
\end{eqnarray*}
and,
$$
x \bar{\star}y := x \ \Sigma_1 \ y + x \ \Sigma_2 \ y. 
$$
We present our 27 relations in the following $9 \times 3$ matrix $[O_{ij}]_{(i:=1, \ldots 9; \ j:=1, \ldots 3)}$. This matrix has 9 blocs of three rows denoted by $(X_{ij})_{(i,j:=1, \ldots, 3)}$.
For instance, $X_{11}$ is the $3 \times 1$-matrix $([O]_{11}, \ [O]_{21}, \ [O]_{31})$, $X_{12}$ is the $3 \times 1$-matrix $([O]_{12}, \ [O]_{22}, \ [O]_{32})$, and so forth. 
\begin{small}
\begin{center}
$
\begin{array}{cccccccccc}
(x \nwarrow_1 y)\nwarrow_1 z &=& x \nwarrow_1(y \ \bar{\star} \ z);& (x \nearrow_1 y)\nwarrow_1 z &=& x \nearrow_1(y \ll z); & (x \ \wedge_1 \ y)\nearrow_1 z &=& x \nearrow_1(y \gg z);  \\
(x \swarrow_1 y)\nwarrow_1 z &=& x \swarrow_1(y \bigwedge z);&
(x \searrow_1 y)\nwarrow_1 z &=& x \searrow_1(y \nwarrow_{12} z);&
(x \ \vee_1 \ y)\nearrow_1 z &=& x \searrow_1(y \nearrow_{12} z); \\
(x \prec_1 y)\swarrow_1 z &=& x \swarrow_1(y \bigvee z); &
(x \succ_1 y)\swarrow_1 z &=& x \searrow_1(y \swarrow_{12} z); &
(x \ \Sigma_1 \ y)\searrow_1 z &=& x \searrow_1 (y \searrow_{12} z);  
\end{array}
$
\end{center}
\begin{center}
$
\begin{array}{cccccccccc}
(x \nwarrow_2 y)\nwarrow_1 z &=& x \nwarrow_2(y \ \Sigma_1 \ z);&
(x \nearrow_2 y)\nwarrow_1 z &=& x \nearrow_2(y \prec_1 z);&
(x \ \wedge_2 \ y)\nearrow_1 z &=& x \nearrow_2(y \succ_1 z);\\
(x \swarrow_2 y)\nwarrow_1 z &=& x \swarrow_2(y \ \wedge_1 \ z);&
(x \searrow_2 y)\nwarrow_1 z &=& x \searrow_2(y \nwarrow_1 z);&
(x \ \vee_2 \ y)\nearrow_1 z &=& x \searrow_2(y \nearrow_1 z);\\
(x \prec_2 y)\swarrow_1 z &=& x \swarrow_2(y \ \vee_1 \ z);&
(x \succ_2 y)\swarrow_1 z &=& x \searrow_2(y \swarrow_1 z);&
(x \ \Sigma_2 \ y)\searrow_1 z &=& x \searrow_2(y \searrow_1 z);
\end{array}
$
\end{center}
\begin{center}
$
\begin{array}{ccccccccc}
(x \nwarrow_{12} y)\nwarrow_2 z &=& x \nwarrow_2(y \ \Sigma_2 \ z);&
(x \nearrow_{12} y)\nwarrow_2 z &=& x \nearrow_2(y \prec_2 z);&
(x \bigwedge y)\nearrow_2 z &=& x \nearrow_2(y \succ_2 z); \\
(x \swarrow_{12} y)\nwarrow_2 z &=& x \swarrow_2(y \ \wedge_2 \ z);&
(x \searrow_{12} y)\nwarrow_2 z &=& x \searrow_2(y \nwarrow_2 z);&
(x \bigvee y)\nearrow_2 z &=& x \searrow_2(y \nearrow_2 z); \\
(x \ll y)\swarrow_2 z &=& x \swarrow_2(y \ \vee_2 \ z);&
(x  \gg  y)\swarrow_2 z &=& x \searrow_2(y \swarrow_2 z);&
(x \ \bar{\star} \ y)\searrow_2 z &=& x \searrow_2(y \searrow_2 z).\\
\end{array}
$
\end{center}
\end{small}
\end{defi}
\begin{theo}
\label{Struct dvh}
Keep notation of Definition \ref{defocto}. If $O$ is an octo-algebra then there exist three quadri-algebras on $O$,
\begin{enumerate}
\item {The depth structure, $(O_d, \ \nearrow_{12}, \ \searrow_{12}, \ \swarrow_{12}, \ \nwarrow_{12})$,}
\item {The vertical structure, $(O_v, \ \prec_1,\ \prec_2,\ \succ_1,\ \succ_2)$,}
\item {The horizontal structure, $(O_h, \ \vee_1,\ \vee_2,\ \wedge_1, \ \wedge _2)$,}
\end{enumerate} 
with: 
\begin{eqnarray*}
x \ \bar{\star} \ y &=& x \nearrow_{12} y + x \searrow_{12} y + x \swarrow_{12} y + x \nwarrow_{12} y, \\ &=& x  \prec_1 y + x \prec_2 y + x \succ_1 y + x \succ_2 y, \\ &=& x  \vee_1  y + x  \vee_2  y  + x  \wedge_1 y + x   \wedge _2 y.
\end{eqnarray*}
The r\^oles of the different binary operations are indicated in the following dictionary. 
If $(Q, \ \nearrow, \ \searrow, \ \swarrow, \ \nwarrow, (\prec, \ \succ, \ \vee, \ \wedge))$ denotes a quadri-algebra obeying thus axioms of Definition \ref{quadri}, then:
\begin{enumerate}
\item {For the quadri-algebra $(O_d, \ \nearrow_{12}, \ \searrow_{12}, \ \swarrow_{12}, \ \nwarrow_{12})$, we get the equivalences: \\ $\nearrow_{12} \ \equiv \ \nearrow, \ \searrow_{12} \ \equiv \ \searrow, \ \swarrow_{12} \ \equiv \ \swarrow, \ \nwarrow_{12} \ \equiv \ \nwarrow$ and $\ll \ \equiv \ \prec, \ \gg \ \equiv \ \succ, \ \bigvee \ \equiv \ \vee  , \ \bigwedge  \ \equiv \ \wedge.$ }
\item {For the quadri-algebra $(O_v, \ \prec_1,\ \prec_2,\ \succ_1,\ \succ_2)$, we have: \\ $\succ_1 \ \equiv \ \nearrow, \ \succ_2 \ \equiv \ \searrow, \ \prec_2 \ \equiv \ \swarrow, \ \prec_1 \ \equiv \ \nwarrow$ and $\ll \ \equiv \ \prec, \ \gg \ \equiv \ \succ, \ \Sigma_2 \ \equiv \ \vee  , \ \Sigma_1  \ \equiv \ \wedge.$}
\item {For the quadri-algebra $(O_h, \ \vee_1,\ \vee_2,\ \wedge_1, \ \wedge _2)$}, we obtain: \\ $\wedge_2 \ \equiv \ \nearrow, \ \vee_2 \ \equiv \ \searrow, \ \vee_1 \ \equiv \ \swarrow, \ \wedge_1 \ \equiv \ \nwarrow$ and $\Sigma_1 \ \equiv \ \prec, \ \Sigma_2 \ \equiv \ \succ, \ \bigvee \ \equiv \ \vee, \ \bigwedge  \ \equiv \ \wedge.$
\end{enumerate} 
\end{theo}
\Proof
Tedious, but no difficulty.
Recall the $9 \times 3$ matrix $[O]$ is a $3 \times 3$ bloc matrix $X$ with each $X_{i,j}$ denoting a bloc composed with three equations. 
\begin{center}
\includegraphics*[width=15cm]{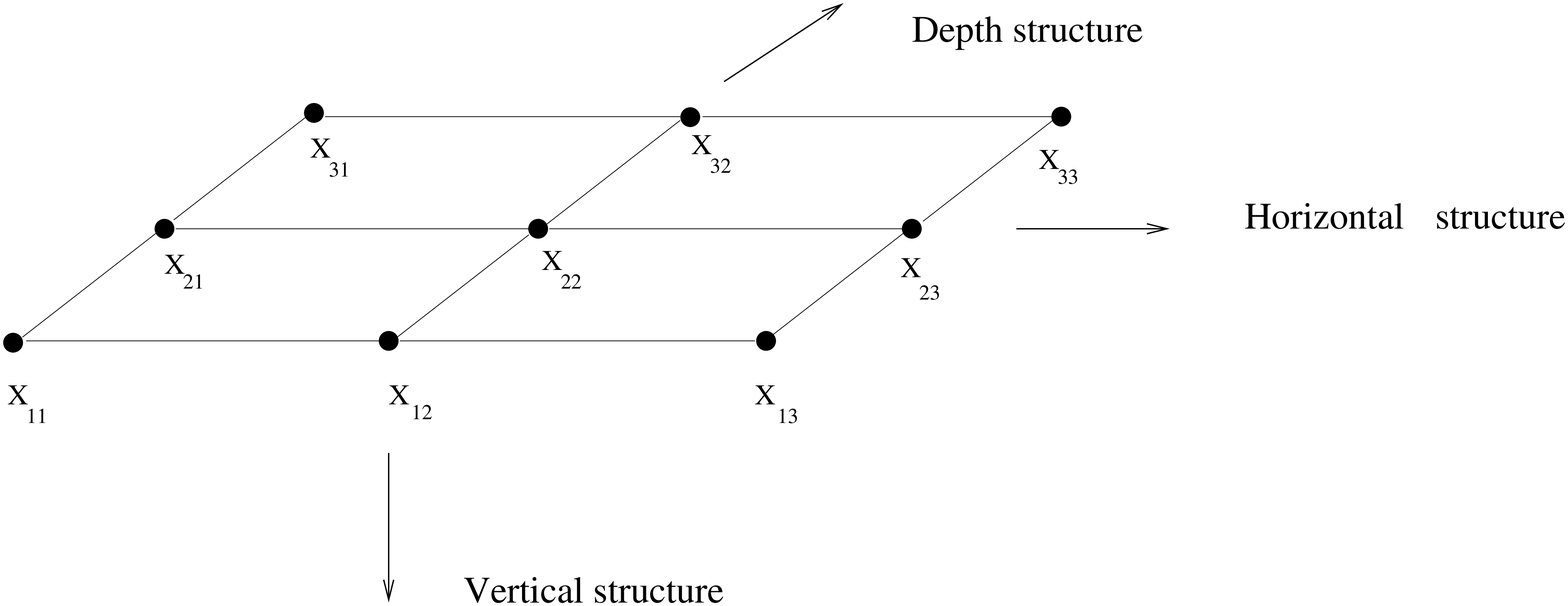}
\end{center}
The depth structure is then obtained by summing as indicated in the above picture. The relations are: $$(O_d)_{ij} := \sum_{k=1}^3 (X_{kj})_i.$$
The vertical structure is obtained by the following relations:
$$(O_v)_{ij} := \sum_{k=1}^3 (X_{ij})_k.$$
The horizontal structure is obtained by the following relations:
$$(O_h)_{ij} := \sum_{k=1}^3 (X_{ik})_j.$$
Up to a permutation of $S_9$, these matrices of relations are the same that the one in Definition \ref{quadri}.
\eproof
\Rk
In terms of categories, Theorem \ref{Struct dvh} implies that the following diagram commutes:
\begin{center}

$
\begin{array}{ccccc}
 \textsf{As.} & \stackrel{f}{\longleftarrow} & \textsf{Quadri.} & \stackrel{f}{\longrightarrow} &  \textsf{As.}\\
& & & & \\ 
\uparrow f & \nwarrow & \uparrow F_d & \nearrow & \uparrow f\\
& & & & \\ 
\textsf{Quadri.} & \stackrel{F_h}{\longleftarrow} & \textsf{Octo.} & \stackrel{F_v}{\longrightarrow} & \textsf{Quadri.} \\ 
 & & & &\\ 
  & \searrow f &  \downarrow   & \swarrow f &  \\
& & & & \\ 
 &   & \textsf{As.} &
\end{array} 
$
\end{center}
\begin{prop} \textbf{[Opposite of an octo-algebra]}  
Let $O(\nearrow_1,\nearrow_2,\searrow_1,\searrow_2,\swarrow_1,
\swarrow_2,\nwarrow_1,\nwarrow_2)$ be an octo-algebra. For all $x,y \in O$, define
the following binary operations as follows:
$$ x \nearrow_{1}^{op} y := y \swarrow_{2} x, \ \ 
  x \nearrow_{2}^{op} y := y \swarrow_{1} x, \ \ 
  x \searrow_{1}^{op} y := y \nwarrow_{2} x, \ \ 
  x \searrow_{2}^{op} y := y \nwarrow_{1} x, \ \ $$ 
 $$ x \swarrow_{1}^{op} y := y \nearrow_{2} x, \ \ 
  x \swarrow_{2}^{op} y := y \nearrow_{1} x, \ \ 
  x \nwarrow_{1}^{op} y := y \searrow_{2} x, \ \ 
  x \nwarrow_{2}^{op} y := y \searrow_{1} x. $$
Then, for all $i=1,2$ we have (the labels being undestood mod 2): 
$$x \prec_{i}^{op} y := y \succ_{i+1} x, \ \ x \succ_{i}^{op} y := y \prec_{i+1} x, \ \ x \wedge_{i}^{op} y := y \vee_{i+1} x, \ \ x \vee_{i}^{op} y := y \wedge_{i+1} x.$$
In addition,
$O(\nearrow_{1}^{op},\nearrow_{2}^{op},\searrow_{1}^{op},\searrow_{2}^{op},\swarrow_{1}^{op},
\swarrow_{2}^{op},\nwarrow_{1}^{op},\nwarrow_{2}^{op})$
is an octo-algebra called the {\it{opposite}} of the octo-algebra $O(\nearrow_{1},\nearrow_{2},\searrow_{1},\searrow_{2},\swarrow_{1},
\swarrow_{2},\nwarrow_{1},\nwarrow_{2})$. 
\end{prop}
\Proof
This comes from the central symmetry of center $[O]_{52}:=(x \searrow_2 y)\nwarrow_1 z= x \searrow_2(y \nwarrow_1 z)$. 
\eproof
\Rk
An octo-algebra coinciding with its opposite is called a {\it{commutative octo-algebra}}. Indeed, for all $x,y$ belonging to a commutative octo-algebra, $x \bar{\star} y = y \bar{\star} x$ holds.

\noindent
There exist 5 other symmetries letting the matrix of relations of an octo-algebra globally invariant. 
\begin{prop} \textbf{[Transpose of an octo-algebra]} 
\label{tr, (xi, (gamma,beta))}
Let $O(\nearrow_1,\nearrow_2,\searrow_1,\searrow_2,\swarrow_1,
\swarrow_2,\nwarrow_1,\nwarrow_2)$ be an octo-algebra. For all $x,y \in O$, define 
the following binary operations as follows:
$$x \nearrow_{1}^t y := x \swarrow_{1} y, \ \ x \nearrow_{2}^t y := x \swarrow_{2} y, \ \ x \swarrow_{1}^t y := x \nearrow_{1} y, \ \ x \swarrow_{2}^t y := x \nearrow_{2} y,$$
and $ x \diamond^t y := x \diamond y$ for $\diamond =
\nearrow_2,\searrow_1,\searrow_2,\swarrow_1,
\swarrow_2,\nwarrow_1,\nwarrow_2.$
Then, for all $i=1,2$, we have: 
$$x \prec_{i}^t y := x \wedge_{i} y, \ \ x \succ_{i}^t y := x \vee_{i} y, \ \ x \wedge_{i}^t y := x \prec_{i} y, \ \ x \vee_{i}^t y := x \succ_{i} y,$$
and,
$$
x \searrow_{12}^t y := x \searrow_{12} y, \ \ x \swarrow_{12}^t y := x \nearrow_{12} y, \ \ x \nearrow_{12}^t y := x \swarrow_{12} y, \ \ x \nwarrow_{12}^t y := x \nwarrow_{12} y. 
$$
In addition,
$O(\nearrow_{1}^t,\nearrow_{2}^t,\searrow_{1}^t,\searrow_{2}^t,\swarrow_{1}^t,
\swarrow_{2}^t,\nwarrow_{1}^t,\nwarrow_{2}^t)$
is an octo-algebra called the {\it{transpose}} of the octo-algebra $O(\nearrow_{1},\nearrow_{2},\searrow_{1},\searrow_{2},\swarrow_{1},
\swarrow_{2},\nwarrow_{1},\nwarrow_{2})$. 
\end{prop}
\Proof
Observe that $[O]_{ij} = [O^t]_{ji}$ on each $3 \times 3$-matrices given by $(X_{11}, \ X_{12}, \ X_{13})$, $(X_{21}, \ X_{22}, \ X_{23})$ and $(X_{31}, \ X_{32}, \ X_{33})$. 
\eproof
\begin{prop} \textbf{[Symmetry $Sym_a$]} 
\label{a, (gamma, (beta,xi))}
Let $O(\nearrow_1,\nearrow_2,\searrow_1,\searrow_2,\swarrow_1,
\swarrow_2,\nwarrow_1,\nwarrow_2)$ be an octo-algebra. For all $x,y \in O$, define 
the following binary operations as follows:
$$x \searrow_{1}^a y := x \nearrow_{2} y, \ \ x \nearrow_{2}^a y := x \searrow_{1} y, \ \ x \swarrow_{1}^a y := x \nwarrow_{2} y, \ \ x \nwarrow_{2}^a y := x \swarrow_{1} y,$$
and let the other ones fixed.
Then, 
$$ x \searrow_{12}^a y  = x \succ_2 y, \ \ x \nearrow_{12}^a y  = x \succ_1 y, \ \ x \swarrow_{12}^a y  = x \prec_2 y, \ \ x \nwarrow_{12}^a y  = x \prec_1 y,$$
and,
$$ x \vee_2 ^a y := x \vee_2 y, \ \ x \vee_1 ^a y := x \wedge_2 y, \ \ x \wedge_1 ^a y := x \wedge_1 y, \ \ x \wedge_2 ^a y := x \vee_1 y.$$
In addition,
$Sym_a(O) := O(\nearrow_{1}^a,\nearrow_{2}^a,\searrow_{1}^a,\searrow_{2}^a,\swarrow_{1}^a,
\swarrow_{2}^a,\nwarrow_{1}^a,\nwarrow_{2}^a)$
is an octo-algebra.
\end{prop}
\Proof
Straightforward.
\eproof
\begin{prop} \textbf{[Symmetry $Sym_b$]} 
\label{b, (gamma, (xi,beta))}
Let $O(\nearrow_1,\nearrow_2,\searrow_1,\searrow_2,\swarrow_1,
\swarrow_2,\nwarrow_1,\nwarrow_2)$ be an octo-algebra. For all $x,y \in O$, define 
the following binary operations as follows:
$$x \searrow_{1}^b y := x \nearrow_{2} y, \ \ x \nearrow_{2}^b y := x \swarrow_{2} y, \ \ x \swarrow_{1}^b y := x \nearrow_{1} y, \ \ x \nwarrow_{2}^b y := x \swarrow_{1} y,$$ $$ \ \ x \nearrow_{1}^b y := x \nwarrow_{2} y, \ \ x \swarrow_{2}^b y := x \searrow_{1} y,$$
and let the two other ones fixed.
Then, 
$$ x \searrow_{12}^b y  = x \succ_2 y, \ \ x \nearrow_{12}^b y  = x \prec_2 y, \ \ x \swarrow_{12}^b y  = x \succ_1 y, \ \ x \nwarrow_{12}^b y  = x \prec_1 y,$$
and,
$$ x \vee_2 ^b y := x \searrow_{12} y, \ \ x \vee_1 ^b y := x \nearrow_{12} y, \ \ x \wedge_1 ^b y := x \nwarrow_{12} y, \ \ x \wedge_2 ^b y := x \swarrow_{12} y,$$
and,
$$ x \prec_2 ^b y := x \vee_1 y, \ \ x \prec_1 ^b y := x \wedge_1 y, \ \ x \succ_1 ^b y := x \wedge_2 y, \ \ x \succ_2 ^b y := x \vee_2 y.$$
In addition,
$Sym_b(O) := O(\nearrow_{1}^b,\nearrow_{2}^b,\searrow_{1}^b,\searrow_{2}^b,
\swarrow_{1}^b, \swarrow_{2}^b,\nwarrow_{1}^b,\nwarrow_{2}^b)$
is an octo-algebra.
\end{prop}
\Proof
Straightforward.
\eproof
\begin{prop} \textbf{[Symmetry $Sym_c$]} 
\label{c, (beta, (gamma,xi))}
Let $O(\nearrow_1,\nearrow_2,\searrow_1,\searrow_2,\swarrow_1,
\swarrow_2,\nwarrow_1,\nwarrow_2)$ be an octo-algebra. For all $x,y \in O$, define 
the following binary operations as follows:
$$x \searrow_{1}^c y := x \swarrow_{2} y, \ \ x \nearrow_{2}^c y := x \searrow_{1} y, \ \ x \swarrow_{1}^c y := x \nwarrow_{2} y, \ \ x \nwarrow_{2}^c y := x \nearrow_{1} y,$$ $$ \ \ x \nearrow_{1}^c y := x \swarrow_{1} y, \ \ x \swarrow_{2}^c y := x \nearrow_{2} y,$$
and let the two other ones fixed.
Then, 
$$ x \searrow_{12}^c y  = x \vee_2 y, \ \ x \nearrow_{12}^c y  = x \vee_1 y, \ \ x \swarrow_{12}^c y  = x \wedge_2 y, \ \ x \nwarrow_{12}^c y  = x \wedge_1 y,$$
and,
$$ x \vee_2 ^c y := x \succ_{2} y, \ \ x \vee_1 ^c y := x \prec_{2} y, \ \ x \wedge_1 ^c y := x \prec_{1} y, \ \ x \wedge_2 ^c y := x \succ_{1} y,$$
and,
$$ x \prec_2 ^c y := x \nearrow_{12} y, \ \ x \prec_1 ^c y := x \nwarrow_{12} y, \ \ x \succ_1 ^c y := x \swarrow_{12} y, \ \ x \succ_2 ^c y := x \searrow_{12} y.$$
In addition,
$Sym_c(O) := O(\nearrow_{1}^c,\nearrow_{2}^c,\searrow_{1}^c,\searrow_{2}^c,
\swarrow_{1}^c, \swarrow_{2}^c,\nwarrow_{1}^c,\nwarrow_{2}^c)$
is an octo-algebra.
\end{prop}
\Proof
Straightforward.
\eproof
\begin{prop} \textbf{[Symmetry $Sym_d$]} 
\label{d, (beta, (xi,gamma))}
Let $O(\nearrow_1,\nearrow_2,\searrow_1,\searrow_2,\swarrow_1,
\swarrow_2,\nwarrow_1,\nwarrow_2)$ be an octo-algebra. For all $x,y \in O$, define 
the following binary operations as follows:
$$x \searrow_{1}^d y := x \swarrow_{2} y, \ \ x \nearrow_{1}^d y := x \nwarrow_{2} y, \ \ x \swarrow_{2}^d y := x \searrow_{1} y, \ \ x \nwarrow_{2}^d y := x \nearrow_{1} y,$$
and let the other ones fixed.
Then, 
$$ x \searrow_{12}^d y  = x \vee_2 y, \ \ x \nearrow_{12}^d y  = x \wedge_2 y, \ \ x \swarrow_{12}^d y  = x \vee_1 y, \ \ x \nwarrow_{12}^d y  = x \wedge_1 y,$$
and,
$$ x \vee_2 ^d y := x \searrow_{12} y, \ \ x \vee_1 ^d y := x \swarrow_{12} y, \ \ x \wedge_1 ^d y := x \nwarrow_{12} y, \ \ x \wedge_2 ^d y := x \nearrow_{12} y,$$
and,
$$ x \prec_2 ^d y := x \succ_{1} y, \ \ x \prec_1 ^d y := x \prec_{1} y, \ \ x \succ_1 ^d y := x \prec_{2} y, \ \ x \succ_2 ^d y := x \succ_{2} y.$$
In addition,
$Sym_d(O) := O(\nearrow_{1}^d,\nearrow_{2}^d,\searrow_{1}^d,\searrow_{2}^d,
\swarrow_{1}^d, \swarrow_{2}^d,\nwarrow_{1}^d,\nwarrow_{2}^d)$
is an octo-algebra.
\end{prop}
\Proof
Straightforward.
\eproof
\subsection{The augmented free octo-algebra and connected Hopf algebra}
Let $V$ be a $k$-vector space. The {\it{free octo-algebra}} $\mathcal{O}(V)$ on $V$ is by definition
an octo-algebra equipped with a map $i: \ V \mapsto \mathcal{O}(V)$ which satisfies the following universal property:
for any linear map $f: V \xrightarrow{} A$, where $A$ is an octo-algebra, there exists a unique octo-algebra morphism $\bar{f}: \mathcal{O}(V) \xrightarrow{} A$ such that $\bar{f} \circ i=f$.

The operad associated with octo-algebras is binary, quadratic and non-symmetric, thus of the form:
$$ \mathcal{O}(V) := \bigoplus_{n \geq 1} \ \mathcal{O}_n \otimes V^{\otimes n} .$$
In particular, the free octo-algebra on one generator $x$ is $ \mathcal{O}(k) := \bigoplus_{n \geq 1} \mathcal{O}_n,$ where $\mathcal{O}_1:= kx$, $\mathcal{O}_2:= k(x \nwarrow_1 x) \oplus k(x \nwarrow_2 x) \oplus k(x \swarrow_1 x) \oplus k(x \swarrow_2 x) \oplus k(x \nearrow_1 x) \oplus k(x \nearrow_2 x) \oplus k(x \searrow_1 x) \oplus k(x \searrow_2 x).$
The space of three variables made out of eight operations is of dimension $2 \times 8^2= 128$. As we have $27$ relations the $k$-vector space $\mathcal{E}_3$ has a dimension equal to $128-27=101.$
Therefore, the sequence associated with the dimensions of $(\mathcal{E}_n)_{n \in \mathbb{N}}$ starts with $1, \ 8, \ 101 \ldots$
Finding the free octo-algebra on one generator is an open problem.
However, a little bit more can be said on the augmented free octo-algebra.

\noindent
To be as self-contained as possible, we introduce some notation to expose a theorem du to Loday \cite{Lodayscd}. 
Recall that a bialgebra $(H, \mu, \ \Delta, \ \eta, \ \kappa)$
is a unital associative algebra $(H, \mu, \ \eta)$ together with co-unital coassocitive coalgebra $(H, \ \Delta, \ \kappa)$. Moreover, it is required that the coproduct $\Delta$
and the counit $\kappa$ are morphisms of unital algebras. A bialgebra is connected if there exists a filtration $(F_rH)_r$ such that $H = \bigcup_r F_rH$, where $F_0H := k1_H$ and for all $r$,
$$ F_rH := \{x \in H; \ \Delta(x) -1_H \otimes x - x\otimes 1_H \in F_{r-1}H \otimes F_{r-1}H \}.$$
Such a bialgebra admits an antipode. Consequently, connected bialgebras are connected Hopf algebras.  

\noindent
Let $P$ be a binary quadratic operad. By a {\it{unit action}} \cite{Lodayscd}, we mean the choice of two linear applications:
$$\upsilon: P(2) \xrightarrow{} P(1), \ \ \ \ \ \varpi:P(2) \xrightarrow{} P(1),$$
giving sense, when possible, to $x \circ 1$ and $1 \circ x$, for all operations $\circ \in P(2)$ and for all $x$ in the $P$-algebra $A$, i.e.,
$x \circ 1 = \upsilon(\circ)(x)$ and $1 \circ x= \varpi(\circ)(x)$.
If $P(2)$ contains an associative operation, say $\star$, then we will require that $x \star 1 := x := 1 \star x$, i.e., $\upsilon(\star) := Id := \varpi(\star)$.
We say that the {\it{unit action}} or the couple $(\upsilon,\varpi)$ is {\it{compatible}} with the relations of the $P$-algebra $A$ if they still hold on $A_+:= k1 \oplus A$.
Let $A$, $B$ be two $P$-algebras. Using the couple $(\upsilon,\varpi)$, we extend binary operations $\circ \in P(2)$ to the $k$-vector space $A \otimes 1.k \oplus k.1 \otimes B \oplus A \otimes B$ by requiring:
\begin{eqnarray}
(a \otimes b) \circ (a' \otimes b') &:= &(a \star a') \otimes (b \circ b') \ \ \ \textrm{if} \ \ b \otimes b' \not= 1 \otimes 1, \\
(a \otimes 1) \circ (a' \otimes 1) &:= &(a \circ a') \otimes 1, \ \ \ \textrm{otherwise}.
\end{eqnarray}
The couple $(\upsilon,\varpi)$ is said to be {\it{coherent}} with the relations of $P$ if $A \otimes 1.k \oplus k.1 \otimes B \oplus A \otimes B$, equipped with these operations is still a $P$-algebra. Observe that a necessary condition for having coherence is compatibility.

\noindent
One of the main interest of these two concepts is the construction of a connected Hopf algebra on the augmented free $P$-algebra. 
\begin{theo} [\textbf{\cite{Lodayscd}}]
\label{Lodaych}
Let $P$ be a binary quadratic operad. Suppose there exists an associative operation in $P(2)$. Any unit action coherent with the relations of $P$ equips the augmented free $P$-algebra $P(V)_+$ on a $k$-vector space $V$ with a
coassociative coproduct $\Delta: P(V)_+ \xrightarrow{}  P(V)_+ \otimes P(V)_+,$
which is a $P$-algebra morphism.
Moreover,
$P(V)_+$ is a connected Hopf algebra.
\end{theo}
\Proof
See \cite{Lodayscd} for the proof. However, we reproduce it to make things clearer.

\noindent
Let $V$ be a $k$-vector space and $P(V)$ be the free $P$-algebra on $V$.
Thanks to the coherence of the unit action, $P(V)_+ \otimes P(V)_+$ is a $P$-algebra. Consider the linear map $\delta: V \xrightarrow{} P(V)_+ \otimes P(V)_+$, given by $v \mapsto 1 \otimes v + v \otimes 1$. Since, $P(V)$ is the free $P$-algebra on $V$, there exists a unique extension of $\delta$ to a morphism of augmented $P$-algebra $\Delta: P(V)_+ \xrightarrow{}  P(V)_+ \otimes P(V)_+$. Now, $\Delta$ is coassociative since the morphisms $(\Delta \otimes id)\Delta$ and 
$(id \otimes \Delta)\Delta$ extend the linear map $ V \xrightarrow{} P(V)_+ ^{\otimes 3}$ which maps $v$ to $1 \otimes 1 \otimes v + 1 \otimes v \otimes 1 + v \otimes 1 \otimes 1$. By unicity of the extension, the coproduct $\Delta$ is coassociative. The bialgebra so obtained is connected. Indeed, by definition, the free $P$-algebra $P(V)$ can be written such as
$P(V) := \oplus_{n \geq 1} \ P(V)_n$, where $P(V)_n$ is the $k$-vector space of products of $n$ elements of $V$. Moreover, we have $\Delta(x \circ y ) := 1 \otimes (x \circ y )  + (x \circ y) \otimes 1 + x \otimes (1 \circ y) + y \otimes (x \circ 1)$, for all $x,y \in P(V)$ and $\circ \in P(2)$. The filtration of $P(V)_+$ is then $FrP(V)_+ = k.1 \oplus \oplus_{1 \leq n \leq r} P(V)_n$.
Therefore, $P(V)_+:=\cup_r \ FrP(V)_+$ and $P(V)_+$ is a connected bialgebra.
\eproof

\noindent
We will use this theorem to show that there exists a connected Hopf algebra on the augmented free octo-algebra as well as on the augmented free commutative octo-algebra.
\begin{prop}
Let $\mathcal{O}(V)$ be the free octo-algebra on a $k$-vector space $V$. Extend the binary operations $\nwarrow_1$ and $\searrow_2$ to $\mathcal{O}(V)_+$ as follows:
$$ x \nwarrow_1 1 := x, \  \ 1 \nwarrow_1 x := 0, \ \ 1 \searrow_2 x := x, \ \ x \searrow_2 1 := 0, \ \ \forall x \in \mathcal{O}(V).$$
In addition, for all operations $\diamond \in \{\searrow_1, \ \nwarrow_2, \ \nearrow_1, \nearrow_2, \swarrow_1, \swarrow_2\}$ choose: 
$ x \diamond 1 := 0 = 1 \diamond x, \ \ \forall x \in \mathcal{O}(V).$ Then,
$$ x \prec_1 1 =x, \ \ 1 \succ_2 x =x, \ \ 1 \vee_1 x := x , \ \ x \wedge_2 1 := x, \ \ x \bar{\star} 1 = x = 1 \bar{\star} x, \ \forall x. \in \mathcal{O}(V)$$ 
Moreover, this choice is coherent.
\end{prop}
\Proof
Observe that we cannot extend $\nwarrow_1$ and $\searrow_2$ to the field $k$, i.e., $1 \nwarrow_1 1$ and $1 \searrow_2 1$ have no sense.
Keep notation introduced in that Subsection.
Firstly, we claim that this choice is compatible. Let $x,y,z \in \mathcal{O}(V)_+$. We have to show for instance that the relation $(x \nwarrow_1 y)\nwarrow_1 z = x \nwarrow_1 (y \bar{\star} z)$ holds in $\mathcal{O}(V)_+$. Indeed, for $x=1$ we get $0=0$. For $y =1$ we get $ x \nwarrow z = x \nwarrow z$ and for $z=1$ we get $ x \nwarrow y = x \nwarrow y$. We do the same thing with the 27 others and quickly found that the augmented free octo-algebra $\mathcal{O}(V)_+$ is still an octo-algebra.

\noindent
Secondly, let us show that this choice is coherent. Let $x_1,x_2, x_3, y_1, y_2, y_3 \in \mathcal{O}(V)_+$. We have to show that, for instance:
$$(Eq. \ [O]_{11}) \ \ \  ((x_1 \otimes y_1) \nwarrow_1 (x_2 \otimes y_2)) \nwarrow_1 (x_3 \otimes y_3) = (x_1 \otimes y_1) \nwarrow_1 ((x_2 \otimes y_2) \bar{\star} (x_3 \otimes y_3)).$$
Indeed, if there exists a unique $y_i =1$, the other belonging to $\mathcal{O}(V)$, then, by definition we get:
$$ (x_1 \bar{\star} x_2 \bar{\star} x_3 ) \otimes ((y_1 \nwarrow_1 y_2) \nwarrow_1  y_3) = (x_1 \bar{\star} x_2 \bar{\star} x_3 ) \otimes (y_1 \nwarrow_1 ( y_2 \bar{\star}  y_3)),$$
which always hold since our choice of the unit action is compatible. Similarly if $y_1=y_2=y_3=1$.
For $y_1=1=y_2$ and $y_3 \in \mathcal{O}(V)$, we get: $0=0$, similarly for $y_1=1=y_3$ and $y_2 \in \mathcal{O}(V)$. If $y_1 \in \mathcal{O}(V)$ and $y_2=1=y_3$, the two hand sides of $(Eq. \ [O]_{11})$ are equal to $(x_1 \bar{\star} x_2 \bar{\star} x_3 ) \otimes y_1$. As a conclusion, $ (Eq. \ [O]_{11})$ holds in $\mathcal{O}(V)  \otimes 1.k \oplus k.1 \otimes \mathcal{O}(V)  \oplus \mathcal{O}(V)  \otimes \mathcal{O}(V) $. Checking the same thing with the 27 others shows that our choice of the unit action is coherent.
\eproof
\begin{coro}
There exists a connected Hopf algebra structure on the augmented free octo-algebra as well as on the augmented free commutative 
octo-algebra.
\end{coro}
\Proof
The first claim comes from the fact that our choice of the unit action is coherent and from Theorem \ref{Lodaych}. For the second remark, observe that our choice is in agreement with the symmetry relations defining a commutative octo-algebra since for instance $x \nwarrow_1 ^{op} 1 := 1 \searrow_2 x :=x$ and $1 \searrow_2 ^{op} x := x \nwarrow_1 1 :=x,$ for all $x \in \mathcal{O}(V)$.
\eproof
\subsection{Construction of octo-algebras via three pairwise commuting Baxter operators}
\noindent
To produce such exotic associative algebras, it suffices to construct a Baxter operator on a quadri-algebra $(Q, \ \nearrow, \ \searrow, \ \swarrow, \ \nwarrow)$. By definition,  
such an operator is a linear map $\xi: Q^{\otimes 2} \xrightarrow{} Q$ such that
for all binary operations $\circ \in \{ \  \nearrow, \ \searrow, \ \swarrow, \ \nwarrow \}$ and all $x,y \in Q$,
$$ \xi(x) \circ \xi(y) := \xi(\xi(x) \circ y + x \circ \xi(y)),$$
hold. In particular, if $\star \longrightarrow \ \nearrow + \searrow + \swarrow + \nwarrow$ denotes the associative operation of $Q$, then
$\xi(x) \star \xi(y) := \xi(\xi(x) \star y + x \star \xi(y)).$
\begin{prop}
Let $(Q, \ \nearrow, \ \searrow, \ \swarrow, \ \nwarrow)$ be a quadri-algebra and $\xi: Q^{\otimes 2} \xrightarrow{} Q$ be a Baxter operator. Define new binary operations as follows:
$x\searrow_1 y := x \searrow \xi(y), \ \
x\searrow_2 y := \xi(x) \searrow y,\ \
x\nearrow_1 y := x \nearrow \xi(y),\ \
x\nearrow_2  y:= \xi(x) \nearrow y,\ \
x\nwarrow_1 y := x \nwarrow \xi(y),\ \
x\nwarrow_2 y := \xi(x) \nwarrow y,\ \
x\swarrow_1  y:= x \swarrow \xi(y),\ \
x\swarrow_2  y:= \xi(x) \swarrow y, $
for all $x,y \in Q$.

\noindent
Then, $(Q, \ \swarrow_1, \ \swarrow_2, \ \nwarrow_1, \ \nwarrow_2, \ \searrow_1, \ \searrow_2, \ \nearrow_1, \ \nearrow_2)$ is an octo-algebra.
\end{prop}
\Proof
Straightforward. Denotes the associative operation of an octo-algebra by $\bar{\star} \longrightarrow \ \swarrow_1 + \swarrow_2 + \nwarrow_1 + \nwarrow_2 + \searrow_1 + \searrow_2 + \nearrow_1 + \nearrow_2$. Fix $x,y,z \in Q$. For instance,
\begin{eqnarray*}
(x \nwarrow_1 y) \nwarrow_1 z &:=& (x \nwarrow \xi(y)) \nwarrow \xi(z),  \\
&:=& x \nwarrow (\xi(y) \star \xi(z)), \\
&:=& x \nwarrow \xi(y \bar{\star} z), \\
&:=& x \nwarrow_1 (y \bar{\star} z).
\end{eqnarray*}
\eproof
\begin{prop}
\label{Prop. 6.10}
Let $\beta, \ \gamma$ and $\xi$ be three pairwise commuting Baxter operators on an associative algebra $(A, \ \mu)$. Then, $\xi$ is a Baxter operator on the 
quadri-algebra $(Q, \ \swarrow, \ \searrow, \ \nearrow, \ \nwarrow)$ where binary operations are given by:
$$x \searrow y := \beta\gamma(x)y, \ \
x \nearrow y := \beta(x)\gamma(y),\ \
x \swarrow y := \gamma(x)\beta(y),\ \
x \nwarrow y := x\beta\gamma(y).$$
\end{prop}
\Proof
Fix $x,y \in Q$. For instance,
$ \xi(x) \searrow \xi(y) := \beta\gamma(\xi(x)) \xi(y) = \xi(\beta\gamma(x)) \xi(y) = \xi(\beta\gamma(x) \xi(y) + \xi(\beta\gamma(x)) y ) = \xi(\beta\gamma(x) \xi(y) + \beta\gamma(\xi(x)) y ) = \xi(x \searrow \xi(y) + \xi(x) \searrow y ).$
\eproof
\begin{coro}
Let $(A, \mu, \Delta_1 \curvearrowright \Delta_2)$ be a $L$-circular $\epsilon$-bialgebra, i.e., both $(A, \mu, \Delta_1)$ and $(A, \mu, \Delta_2)$ are $\epsilon$-bialgebras and,
$$ (id \otimes \Delta_1)\Delta_2=(\Delta_1 \otimes id)\Delta_2=(id \otimes\Delta_2)\Delta_1.$$
Then, the three Baxter operators $\beta_1, \ \gamma_2, \ \beta_2: \textsf{End}(A) \xrightarrow{} \textsf{End}(A)$ commute pairwise and the linear map $\gamma_2$ is a Baxter operator on the quadri-algebra $\textsf{End}(A)$ whose binary operations are given by:
$$T \searrow S := \beta_1\gamma_2(T)S, \ \
T \nearrow S := \beta_1(T)\gamma_2(S),\ \
T \swarrow S := \gamma_2(T)\beta_1(S),\ \
T \nwarrow S := T\beta_1\gamma_2(S).$$
\end{coro}
\Proof
Since $\Delta_2$ is coassociative, $\beta_2$ and $\gamma_2$ commute. Since $(\Delta_1 \otimes id)\Delta_2=(id \otimes\Delta_2)\Delta_1$, $\beta_1$ and $\gamma_2$ commute. Since $(id \otimes \Delta_1)\Delta_2=(id \otimes\Delta_2)\Delta_1$, $\beta_1$ and $\beta_2$ commute. The corollary holds by applying
Proposition \ref{Prop. 6.10}.
\eproof
\Rk
In additon to the entanglement condition between $\Delta_1$ and $\Delta_2$, observe the importance of the circular condition for having three pairwise commuting Baxter operators.
\Rk
We started with the Baxter operator $\xi$ acting on the quadri-algebra
defined by the pair $(\beta, \gamma)$ and found eight operations defining an octo-algebra, say $O$. If the Baxter operator $\xi$ acts on the quadri-algebra
defined by the pair $(\gamma, \beta)$, we will find the transpose of the octo-algebra $O$, see Proposition \ref{tr, (xi, (gamma,beta))}. Now, if the Baxter operator $\gamma$ acts on the quadri-algebra
defined by the pair $(\beta, \xi)$, then we will find the octo-algebra of Proposition \ref{a, (gamma, (beta,xi))}, i.e. $Sym(O)_a$. If it acts on the pair $(\xi, \beta)$, then we will found the octo-algebra of Proposition \ref{b, (gamma, (xi,beta))}, i.e. $Sym(O)_b$. If the Baxter operator $\beta$ acts on the quadri-algebra
defined by the pair $(\gamma, \xi)$, then we will find the octo-algebra of Proposition \ref{c, (beta, (gamma,xi))}, i.e. $Sym(O)_c$. If it acts on the pair $(\xi, \gamma)$, we will find the octo-algebra of Proposition \ref{d, (beta, (xi,gamma))}, i.e. $Sym(O)_d$. 

\Rk
For the time being, except the trivial ones, i.e., $\Delta_2$ proportional to $\Delta_1$, we do not know any examples of $L$-circular bialgebras.
\subsection{Another construction}
\begin{prop}
\label{treeocto}
Let $(E, \ \prec_{e}, \ \succ_{e}), \ (F, \ \prec_{f}, \ \succ_{f})$ and  $(G, \ \prec_{g}, \ \succ_{g})$ be three dendriform dialgebras.
For all $x,x' \in E$, $y,y' \in F$, $z,z' \in G$, consider the following operations on $E \otimes F \otimes G$:
\begin{eqnarray*}
(x \otimes y \otimes z) \nwarrow_{1} (x' \otimes y' \otimes z') &:=& (x \prec_{e} x') \otimes (y \prec_{f} y') \otimes (z \prec_{g} z'), \\
(x \otimes y \otimes z) \nwarrow_{2} (x' \otimes y' \otimes z') &:=& ( x \prec_{e} x') \otimes (y \prec_{f} y') \otimes (z \succ_{g} z'), \\
(x \otimes y \otimes z) \searrow_{1} (x' \otimes y' \otimes z') &:=& ( x \succ_{e} x') \otimes (y \succ_{f} y') \otimes (z \prec_{g} z'), \\
(x \otimes y \otimes z) \searrow_{2} (x' \otimes y' \otimes z') &:=& ( x \succ_{e} x') \otimes (y \succ_{f} y') \otimes (z \succ_{g} z'), \\
(x \otimes y \otimes z) \nearrow_{1} (x' \otimes y' \otimes z') &:=& ( x \succ_{e} x') \otimes (y \prec_{f} y') \otimes (z \prec_{g} z'), \\
(x \otimes y \otimes z) \nearrow_{2} (x' \otimes y' \otimes z') &:=& ( x \succ_{e} x') \otimes (y \prec_{f} y') \otimes (z \succ_{g} z'), \\
(x \otimes y \otimes z) \swarrow_{1} (x' \otimes y' \otimes z') &:=& ( x \prec_{e} x') \otimes (y \succ_{f} y') \otimes (z \prec_{g} z'), \\
(x \otimes y \otimes z) \swarrow_{2} (x' \otimes y' \otimes z')&:=& (x \prec_{e} x') \otimes (y \succ_{f} y') \otimes (z \succ_{g} z').
\end{eqnarray*}
Equipped with these eight operations, the $k$-vector space $E \otimes F \otimes G$ is an octo-algebra.
The depth structure is described by the following four operations:
\begin{eqnarray*}
(x \otimes y \otimes z) \nwarrow_{12} (x' \otimes y' \otimes z') &:=& (x \prec_{e} x') \otimes (y \prec_{f} y') \otimes (z \star_{g} z'), \\
(x \otimes y \otimes z) \searrow_{12} (x' \otimes y' \otimes z') &:=& ( x \succ_{e} x') \otimes (y \succ_{f} y') \otimes (z \star_{g} z') \\
(x \otimes y \otimes z) \nearrow_{12} (x' \otimes y' \otimes z') &:=& ( x \succ_{e} x') \otimes (y \prec_{f} y') \otimes (z \star_{g} z'), \\
(x \otimes y \otimes z) \swarrow_{12} (x' \otimes y' \otimes z') &:=& ( x \prec_{e} x') \otimes (y \succ_{f} y') \otimes (z \star_{g} z').
\end{eqnarray*}
The vertical structure is described by the following four operations:
\begin{eqnarray*}
(x \otimes y \otimes z) \succ_{1} (x' \otimes y' \otimes z') &:=& (x \succ_{e} x') \otimes (y \star_{f} y') \otimes (z \prec_{g} z'), \\
(x \otimes y \otimes z) \succ_{2} (x' \otimes y' \otimes z') &:=& ( x \succ_{e} x') \otimes (y \star_{f} y') \otimes (z \succ_{g} z') \\
(x \otimes y \otimes z) \prec_{1} (x' \otimes y' \otimes z') &:=& ( x \prec_{e} x') \otimes (y \star_{f} y') \otimes (z \prec_{g} z'), \\
(x \otimes y \otimes z) \prec_{2} (x' \otimes y' \otimes z') &:=& ( x \prec_{e} x') \otimes (y \star_{f} y') \otimes (z \succ_{g} z').
\end{eqnarray*}
The horizontal structure is described by the following four operations:
\begin{eqnarray*}
(x \otimes y \otimes z) \wedge_{1} (x' \otimes y' \otimes z') &:=& (x \ \star_{e} x') \otimes (y \prec_{f} y') \otimes (z \prec_{g} z'), \\
(x \otimes y \otimes z) \wedge_{2} (x' \otimes y' \otimes z') &:=& ( x \star_{e} x') \otimes (y \prec_{f} y') \otimes (z \succ_{g} z') \\
(x \otimes y \otimes z) \vee_{1} (x' \otimes y' \otimes z') &:=& ( x \star_{e} x') \otimes (y \succ_{f} y') \otimes (z \prec_{g} z'), \\
(x \otimes y \otimes z) \vee_{2} (x' \otimes y' \otimes z') &:=& ( x \star_{e} x') \otimes (y \succ_{f} y') \otimes (z \succ_{g} z').
\end{eqnarray*}
The associative structure is then:
$$
(x \otimes y \otimes z) \bar{\star} (x' \otimes y' \otimes z') := ( x \star_{e} x') \otimes (y \star_{f} y') \otimes (z \star_{g} z'),$$
as expected.
\end{prop}
\Proof
Straightforward.
\eproof
\begin{exam}{}
Denote by $T$, the dendriform dialgebra of planar (rooted) binary trees. Consider three copies of $T$, say $T_1$, $T_2$ and $T_3$. Graphically such structures can be represented by placing planar binary trees on a marked directed triangle. On a vertex $i:=1,2,3$, a tree from $T_i$ is placed. The operations of the octo-algebra $T_1 \otimes T_2 \otimes T_3$ are then well defined.
\end{exam}
\begin{center}
\includegraphics*[width=5cm]{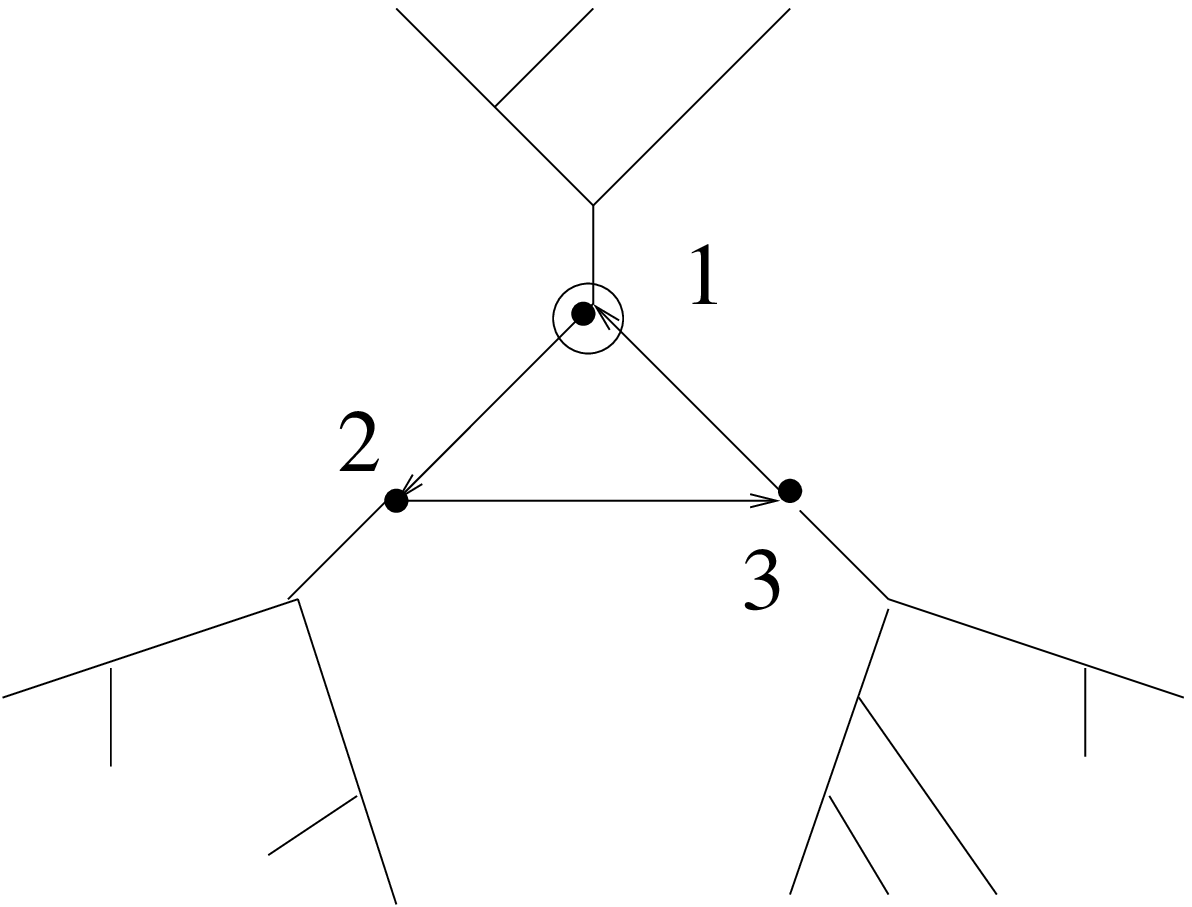}

\begin{scriptsize}
\textbf{Example of planar binary trees on a marked directed triangle.}
\end{scriptsize}
\end{center}
\section{Leibniz algebras and Poisson dialgebras from right Baxter operators }
\subsection{Theoretical constructions}
We end this work on Baxter operators by showing another application to construct particular Leibniz algebras and Poisson dialgebras. Leibniz algebras have been introduced by Loday in \cite{LodayLeib} to understand periodicity phenomena in $K$-theory, see also \cite{Loday}. They are anti-symmetric version of Lie algebras. Similarly, Poisson algebras are the data of a commutative algebra and a Lie bracket verifying compatibility relations. Poisson dialgebras, introduced in \cite{LodayRonco} are the data of an associative algebra, not necessarily commutative, and a Leibniz bracket verifying compatibility conditions, see also \cite{dipt}. A way to obtain Leibniz algebras is to start with associative dialgebras, $k$-vector spaces $D$ equipped with two associative products,
$\vdash$ and $\dashv$, such that for all $x,y,z \in D$, 
\begin{enumerate}
\item{$x \dashv (y \dashv z) = x \dashv (y \vdash z),$}
\item{$(x \vdash y) \dashv z = x \vdash (y \dashv z),$}
\item{$(x \vdash y) \vdash z = (x \dashv y) \vdash z.$}
\end{enumerate} 
Associative dialgebras are then Leibniz algebras by defining the bracket
$[x,y ] := x \dashv y - y \vdash x$, for all $x,y \in D$. This bracket verifies:
$$[[x,y ] ,z ] = [[ x,z ] ,y ] + [ x,[ y  ,z ]].$$
A way to obtain Poisson dialgebras is to start with associative trialgebras. An associative trialgebra is a $k$-vector space $(T, \ \dashv, \ \vdash, \ \perp)$ such that $(T, \ \dashv, \ \vdash)$ is an
associative dialgebra, $(T, \ \perp)$ an associative algebra and,
\begin{enumerate}
\item{$(x \dashv y) \dashv z = x \dashv (y \perp z),$}
\item{$(x \perp y) \dashv z = x \perp (y \dashv z),$}
\item{$(x \dashv y) \perp z = x \perp (y \vdash z),$}
\item{$(x \vdash y) \perp z = x \vdash (y \perp z),$}
\item{$(x \perp y) \vdash z = x \vdash (y \vdash z).$}
\end{enumerate} 
The operads associated with associative di and trialgebras are Koszul dual respectively to the operads associated with dendriform dialgebras and dendriform trialgebras \cite{LodayRonco, KEF, Lertribax}.
We give now a method to construct such types of algebras via right Baxter operators. For that, we define a {\it{$\epsilon[R]$-dialgebra}} $(A, \ \mu, \ \Delta_\vdash, \ \Delta_\dashv)$ which is a $k$-vector space such that $(A, \ \mu, \ \Delta_\vdash)$ and $(A, \ \mu, \ \Delta_\dashv)$ are both $\epsilon[R]$-bialgebras and $(A, \ \Delta_\vdash, \ \Delta_\dashv)$ is a coassociative codialgebra.
\begin{theo}
Let $(A, \ \mu, \ \Delta_\vdash, \ \Delta_\dashv)$ be a $\epsilon[R]$-dialgebra.  Then, $(\textsf{End}(A), \ \overrightarrow{*}_{\gamma_\vdash}, \ \overrightarrow{*}_{\gamma_\dashv})$ is an associative dialgebra. Equipped with the dissymmetric commutator:
$$ [T,S] := T \overrightarrow{*}_{\gamma_\vdash} S - S \overrightarrow{*}_{\gamma_\dashv} T,$$
$(\textsf{End}(A), \ [\cdot, \ \cdot])$ is a Leibniz algebra.
\end{theo}
\Proof
Fix $R,T, S \in \textsf{End}(A)$ and $a \in A$.
The proof can be done without any tedious computations by observing that for computing 
$\gamma_1(T)\gamma_2(S)(a)$ only the map $(id \otimes \Delta_1)\Delta_2$ is required, but for computing $\gamma_1(T \gamma_2(S))(a)$ only the map $(\Delta_2 \otimes id)\Delta_1$ is needed.
For instance, on the one hand,
\begin{eqnarray*}
R \overrightarrow{*}_{\gamma_\vdash} ( S \overrightarrow{*}_{\gamma_\dashv} T ) &:=& R \gamma_\vdash(S \gamma_\dashv(T)).
\end{eqnarray*}
Apply to an element $a \in A$, this equation will only require $(\Delta_\dashv \otimes id)\Delta_\vdash$.
On the other hand, 
\begin{eqnarray*}
(R \overrightarrow{*}_{\gamma_\vdash}  S) \overrightarrow{*}_{\gamma_\vdash} T  &:=& R \gamma_\vdash(S) \gamma_\vdash(T).
\end{eqnarray*}
Apply to an element $a \in A$, this equation will only require $(id \otimes \Delta_\vdash)\Delta_\vdash$. However, $(\Delta_\dashv \otimes id)\Delta_\vdash = (\Delta_\vdash \otimes id)\Delta_\vdash$ since it is an axiom of coassociative codialgebras. Therefore,
 \begin{eqnarray*}
R \overrightarrow{*}_{\gamma_\vdash} ( S \overrightarrow{*}_{\gamma_\dashv} T ) &=& (R \overrightarrow{*}_{\gamma_\vdash}  S )\overrightarrow{*}_{\gamma_\vdash} T = R \overrightarrow{*}_{\gamma_\vdash} ( S \overrightarrow{*}_{\gamma_\vdash} T ).
\end{eqnarray*}
Doing the same thing with the three other equations shows that $(\textsf{End}(A), \ \overrightarrow{*}_{\gamma_\vdash}, \ \overrightarrow{*}_{\gamma_\dashv})$ is an associative dialgebra, with $\overrightarrow{*}_{\gamma_\vdash}$ playing the r\^ole of $\dashv$ and conversely.
\eproof
\begin{coro}
Let $(A, \ \mu, \ \Delta_\vdash, \ \Delta_\dashv, \ \Delta_\perp)$ be a $\epsilon[R]$-trialgebra. Then, $(\textsf{End}(A), \ \overrightarrow{*}_{\gamma_\vdash}, \ \overrightarrow{*}_{\gamma_\dashv}, \ \overrightarrow{*}_{\gamma_\perp})$ is an associative trialgebra. Equipped with the dissymmetric commutator:
$$ [T,S] := T \overrightarrow{*}_{\gamma_\vdash} S - S \overrightarrow{*}_{\gamma_\dashv} T,$$
$(\textsf{End}(A), \ [\cdot, \ \cdot], \ \overrightarrow{*}_{\gamma_\perp})$ is a Poisson dialgebra.
\end{coro}
\Proof
The proof is similar to the previous theorem.
\eproof
\subsection{Constructions from coassociative coverings of directed graphs}
\begin{prop}
\label{coditri}
Let $V$ be a non-empty set. Suppose there exist three coproducts $\Delta_1$, $\Delta_2$ and $\Delta_3$ on $kV$, the free $k$-vector space spanned by $V$, verifying co-trialgebras axioms. Extend them on the free associative $k$-algebra $As(V)$ generated by $V$ as follows,
$$\Delta_{i, \sharp}(v_1 \ldots v_n) := v_1 \ldots v_{n-1}\Delta_{i}(v_n),$$
for all $i=1,2,3$.
Then, $(As(V), \ \Delta_{1, \sharp}, \ \Delta_{2, \sharp}, \ \Delta_{3, \sharp})$ is a $\epsilon[R]$-trialgebra. 
\end{prop}
\Proof
Straightforward.
\eproof

\noindent
Apply Proposition \ref{coditri} to co-dialgebras and co-trialgebras constructed from the notion of coassociative coverings of directed graphs \cite{dipt}. 

\Rk
What has been done in this section can be also realised with left versions.
\section{Conclusion and opennings}
Some mathematical objects can be described with the help of $k$-vector spaces, say $A$, equipped with co-operations and associative products verifying compatibility conditions. These give birth to $\epsilon$-bialgebras and their left and right versions. The use of coproducts naturally endow the
$k$-vector space $\textsf{End}(A)$ with other operations called convolution products. It turns out that studying the placings of operators belonging to $\textsf{End}(A)$, for instance at the right or the left hand side of such convolution products can be coded through the shifts $\gamma$ and $\beta$ and generate special (associative) algebras of a certain type. Therefore the operations coded by these shifts, (or placing operators $\gamma$ and $\beta$ and their combinaisons $\alpha\beta$, $\beta^2$ and so forth)  
give birth to binary quadratic and non-symmetric operads. For the time being, we do not know neither what the free objects associated with these operads are nor if they are Koszul.

\bibliographystyle{plain}
\bibliography{These}

\begin{thebibliography}{10}

\bibitem{Aguiar}
M.~{\textsc{Aguiar}}.
\newblock Infinitesimal bialgebras, pre-lie and dendriform algebras.
\newblock {\em to appear in ``{H}opf algebras: {P}roceedings from an
  {I}nternational {C}onference held at {D}e{P}aul {U}niversity'';
  arXiv:math.QA/0211074}.

\bibitem{AguiarLoday}
M.~{\textsc{Aguiar}} and J.-L. {\textsc{Loday}}.
\newblock Quadri-algebras.
\newblock {\em arXiv:math.QA/0309171}.

\bibitem{Baxter}
G.~{\textsc{Baxter}}.
\newblock An analytic problem whose solution follows from a simple algebraic
  identity.
\newblock {\em Pacific J. Math.}, 10:731--742, 1960.

\bibitem{KEF}
K.~{\textsc{Ebrahimi-Fard}}.
\newblock Loday-type algebras and the {R}ota-{B}axter relation.
\newblock {\em Letters in Mathematical Physics}, 61(2):139--147, 2002.

\bibitem{Fresse}
B.~{\textsc{Fresse}}.
\newblock Koszul duality of operads and homology of partition posets.
\newblock {\em Preprint 2002}.

\bibitem{GK}
V.~{\textsc{Ginzburg}} and M.~{\textsc{Kapranov}}.
\newblock Koszul duality for operads.
\newblock {\em Duke Math. J.}, 76(1):203--272, 1994.

\bibitem{Rota}
S.A. {\textsc{Joni}} and G.-C. {\textsc{Rota}}.
\newblock Coalgebras and bialgebras in combinatorics.
\newblock {\em Stud. Appl. Math.}, 61:93--139, 1979.

\bibitem{Khovanov}
M.~{\textsc{Khovanov}}.
\newblock A categorification of the {J}ones polynomial.
\newblock {\em eprint arXiv:math.QA/9908171}.

\bibitem{Lertribax}
Ph. {\textsc{Leroux}}.
\newblock Ennea-algebras.
\newblock {\em eprint, arXiv:math.QA/0309213, (version 2)}.

\bibitem{dipt}
Ph. {\textsc{Leroux}}.
\newblock From entangled codipterous coalgebras to coassociative manifolds.
\newblock {\em eprint arXiv:math.QA/0301080}.

\bibitem{perorb1}
Ph. {\textsc{Leroux}}.
\newblock Periodic orbits, coassociative grammar and quantum random walk over
  $\mathbb{Z}$.
\newblock {\em eprint arXiv:quant-ph/ 0209100}.

\bibitem{codialg1}
Ph. {\textsc{Leroux}}.
\newblock Tiling the $(n^2,1)$-{D}e-{B}ruijn graph with $n$ coassociative
  coalgebras.
\newblock {\em eprint arXiv:math.QA/ 0209108}.

\bibitem{Coa}
Ph. {\textsc{Leroux}}.
\newblock An algebraic framework of weighted directed graphs.
\newblock {\em To appear in Int. J. Math. Math. Sci.}, 2003.

\bibitem{Lodayscd}
J.-L. {\textsc{Loday}}.
\newblock Scindement d'associativit\'e et alg\`ebres de {H}opf.

\bibitem{LodayLeib}
J.-L. {\textsc{Loday}}.
\newblock Une version non commutative des alg\`ebres de {L}ie: Les alg\`ebres
  de {L}eibniz.
\newblock {\em L'Enseignement Math.}, 39:269--293, 1993.

\bibitem{Loday}
J.-L. {\textsc{Loday}}.
\newblock Dialgebras.
\newblock {\em in Dialgebras and related operads, Lecture Notes in Math.},
  1763:7--66, 2001.

\bibitem{LodayRonco}
J.-L. {\textsc{Loday}} and M.~{\textsc{Ronco}}.
\newblock Trialgebras and families of polytopes.
\newblock {\em eprint arXiv:math.QA/0205043}.

\bibitem{LR1}
J.-L. {\textsc{Loday}} and M.~{\textsc{Ronco}}.
\newblock Alg\`ebres de {H}opf colibres.
\newblock {\em C.R.Acad.Sci Paris}, 337(Ser. I):153--158, 2003.

\bibitem{Pir}
T.~{\textsc{Pirashvili}}.
\newblock Sets with two associative operations.
\newblock {\em Preprint 2003}.

\bibitem{Richter}
B.~{\textsc{Richter}}.
\newblock Dialgebren, {D}oppelagebren und ihre {H}omologie.
\newblock {\em {D}iplomarbeit, {B}onn {U}niversit\"at, unpublished}.

\bibitem{Rota1}
G.-C. {\textsc{Rota}}.
\newblock Baxter operators, an introduction.
\newblock {\em in {G}ian-{C}arlo {R}ota on {C}ombinatorics: {I}ntroductory
  papers and commentaries ({J}.{P}.{S}. {K}ung, {E}d.) {B}irkhauser, {B}oston},
  pages 504--512, 1995.

\bibitem{Voi}
D.~{\textsc{Voiculescu}}.
\newblock Free analysis question {I}: Duality transform for the coalgebra of
  $\partial_{X:B}$.
\newblock {\em arXiv:math.QA/0306172}.

\end{thebibliography}

\end{document}